\documentclass[10pt]{amsart}

\usepackage{amsmath,amsthm,amsopn,amsfonts,amssymb,color}
\usepackage[colorlinks=true,linkcolor=blue,urlcolor=blue]{hyperref} 
\usepackage{cancel}
 \usepackage{amsmath,mathtools}

\usepackage{graphicx}
\usepackage{float}
\usepackage{amsmath}
\usepackage{fancyhdr}
\usepackage{epstopdf}
\usepackage{amsfonts}
\usepackage{amsmath,amsthm}

\usepackage{booktabs}
\usepackage[USenglish]{babel}
\usepackage{cancel}
\usepackage{times}
\usepackage{mathrsfs}
\usepackage{multirow,multicol}
\usepackage{amsfonts}
\usepackage{amsthm}
\usepackage{amsmath}
\usepackage{amsmath,amssymb}
\newtheorem{teo}{Theorem}[section]
\newtheorem{prop}[teo]{Proposition}
\newtheorem{defin}[teo]{Definition}
\newtheorem{remark}[teo]{Remark}

\newtheorem{example}[teo]{Example}
\newtheorem{examples}[teo]{Examples}

\newcommand{\mres}{\mathbin{\vrule height 1.6ex depth 0pt width
0.13ex\vrule height 0.13ex depth 0pt width 1.3ex}}

\newcommand\rn{\mathbb{R}^{N}}

\newcommand\R{\mathbb{R}}

\DeclareMathOperator*{\essinf}{ess \ inf}

\def\elle#1{L^{#1}(\Omega)}

\def\lio{L^{\infty}(\Omega)}
\def\w{H_0^{1}(\Omega)}
\def\io{\int_{\Omega}}
\def\norma#1#2{\|#1\|_{\lower 4pt \hbox{$\scriptstyle #2$}}}
\def\un{u_n}

\def\D{\nabla}

\def\vp{\varphi}

\def\pn{p_n}

\def\finedim
\def\gw{G_{\tilde{k}}(w_n)}

\def\R{\mathbb{R}}

\def\elle#1{L^{#1}(\Omega)}

\def\w{W_0^{1,2}(\Omega)}

\def\w1{W_0^{1,1}(\Omega)}

\def\eps{\varepsilon}
\def\lio{L^{\infty} (\Omega)}

\def\dys{\displaystyle}

\def\w{H_0^{1}(\Omega)}

\def\be{\begin{equation}}
\def\ee{\end{equation}}
\def\bc{\begin{cases}}
\def\ec{\end{cases}}

\def\L1{\mathbb{X} (\Omega)}
\def\be{\begin{equation}}
\def\ee{\end{equation}}

  \usepackage{color}
  
 \numberwithin{equation}{section}

\usepackage{todonotes}

\title[The best approximation of a given function in $L^2$-norm]{The best approximation 
of a given function in $L^2$-norm by Lipschitz functions with gradient constraint}

\author[S. Buccheri]{Stefano Buccheri}
\address{S. Buccheri. 
Institut de Recherche en Math\'ematique et Physique,  Université catholique de Louvain
Chemin du Cyclotron 2, bte L7.01.01
1348 Louvain-la-Neuve, 
Belgium}
\email{buccheristefano@gmail.com}

\author[T. Leonori]{Tommaso Leonori}
\address{T. Leonori. Dipartimento di Scienze di Base e Applicate per l'Ingegneria, 
Universit\`a di Roma "{}Sapienza"{},
Via Antonio Scarpa 10, 00161 Roma, Italia.}
\email{tommaso.leonori@uniroma1.it}

\author[J. D. Rossi]{Julio D. Rossi}
\address{J. D. Rossi. Departamento de Matem\'atica, FCEyN, Universidad de Buenos Aires,
 Pabell\'on I, Ciudad Universitaria (1428),
Buenos Aires, Argentina.}
\email{jrossi@dm.uba.ar}

\begin{document}

\date{\today}

\begin{abstract} The starting point of this paper is the study of the asymptotic behavior, as $p\to\infty$, of the following minimization problem
\be\label{abstract}
\min\left\{\frac1{p}\io|\nabla v|^{p}+\frac12\io(v-f)^2 \,, \quad \ v\in W^{1,p} (\Omega)\right\}.
\ee
We show that the limit problem provides the best approximation, in the $L^2$-norm, of the datum $f$ among all Lipschitz functions with Lipschitz constant less or equal than one. Moreover such approximation verifies a suitable PDE
in the viscosity sense. 

After the analysis of the model problem \eqref{abstract}, we consider the asymptotic behavior of a related family of nonvariational equations  and, finally, we also deal with some functionals involving the $(N-1)$-Hausdorff measure of the jump set of the function.
\end{abstract}

\maketitle 
%\tableofcontents

\section{Introduction}

Let us assume that $\Omega\subset\mathbb{R}^N$ is a bounded open set with smooth, say $C^1$,  boundary.
The main goal of this work is to study the optimal approximation in the $L^2$-norm of a given function $f\in \elle 2$ by functions in 
$W^{1,\infty} (\Omega)$ with a constraint on the gradient. Namely, we consider the following minimization problem
\be\label{minlim.intro}
\inf_{\mathbb{X}(\Omega) }  J_{\infty}(v) \qquad  \mbox{with} \qquad \ J_{\infty}(v)=\frac12 \io(v-f)^2, 
\ee
where 
the set   $\mathbb{X} (\Omega)$ is given by
$$
\mathbb{X} (\Omega) := \Big\{ v \in W^{1,\infty} (\Omega) \,:\, |Dv | \leq 1 \ \text{ for } \ a.e. \ x \in \Omega \Big\}\,.
$$
In order to obtain a minimizer to \eqref{minlim.intro} one can argue by direct methods
(taking a minimizing sequence), 
 or rather noticing that such a problem appears naturally as $\Gamma$-limit, for $n\to\infty$, of the following classical energy functional
\be\label{neum1}
J_{n}(v)=\frac1{p_n}\io|\nabla v|^{p_n}+\frac12\io(v-f)^2 \,, \ \ \qquad \ v\in W^{1,p_n} (\Omega)\,,
\ee
where $p_n$ is a sequence of numbers that diverges to $+\infty$.
Setting  $u_n$ the minimizer of $J_{n}$, the heuristic of the limiting process is that the measure of the set $\{|\nabla u_n|>1\}$ has to go to zero in order to keep the energy bounded (notice that $J_{n}(u_n)\le J_n(0)=\|f\|_{\elle 2}^2$). Therefore, we expect that $u_n$ converge to the    unique minimizer $u\in\mathbb{X} (\Omega) $ of $J_\infty$ (see Proposition \ref{auxlem} below).  

 An interesting feature of the minimizer of \eqref{minlim.intro}
that we find here, is that it satisfies a kind of representation formula by cones inside the region $\{u\neq f\}$.

Before stating our first result, in order to gain some intuition, let us consider three concrete examples in the special case $\Omega = (-1,1)$.

\begin{itemize}
\item  Let $f_1 (x) = k \chi_{(-r,r)}$ with $r \in (0,1]$ and $k\in \mathbb{R}$. 
Observe that if $r=1$ then $f_1 \in \mathbb{X} (-1,1)$ and so $u_1\equiv f_1$
is the minimizer. Otherwise, it is not hard to see that 
$$
u_1 =  \min \Big\{ (r+\frac{k}2-|x|)_+, k \Big\}\,
$$
is the minimizer. In Figure 1  the case $r>\frac{k}2$ is represented.

\medskip 

\centerline{
\includegraphics[width=0.5\textwidth]{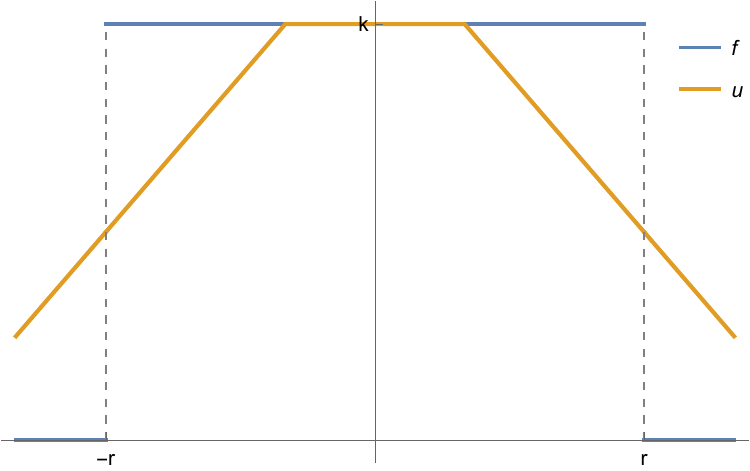} 
}
\centerline{Figure 1.}

\bigskip 

\item Let  $f_2 (x) = 2|x|$, then the minimizer is given by 
$$u_2 (x) =|x|+\frac{1}{2}.$$
 See Figure 2.

\medskip

\centerline{
\includegraphics[width=0.5\textwidth]{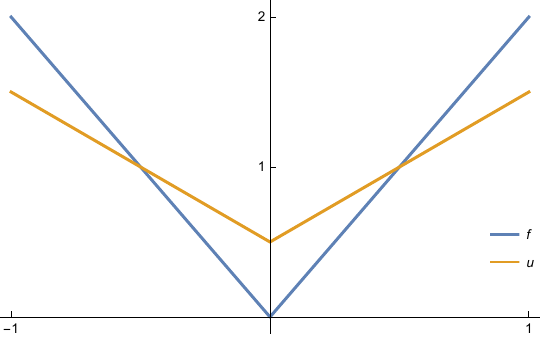}}
\centerline{Figure 2.}

\bigskip

\item Let $f_3(x)=\sqrt{|x|}$, then is is not hard to see that the minimizer is
$$
u_3 (x)=
\begin{cases}
x+ \frac29 & \mbox{for} \ |x|\leq \frac49\,,\\
\sqrt{|x|}& \mbox{for} \ |x|\geq \frac49.
\end{cases}
$$
See Figure 3.

\medskip

\centerline{
\includegraphics[width=0.5\textwidth]{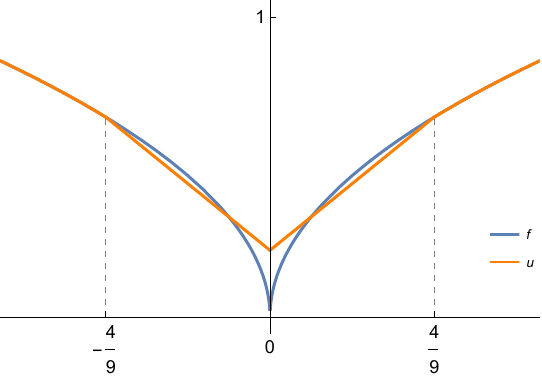}}
 \centerline{Figure 3.}
 
\bigskip

 \end{itemize}

In the three examples above (see Section \ref{exa} for the explicit computations) the behavior of the solution $u_i$, $i=1,2,3$  follows the same idea: any $u_i$ tries to be as close as possible to the datum $f_i$ as long as $f_i$ is smooth and has gradient bounded by one. Otherwise, the best that the solution can do in order to minimize the $L^2$ norm of the difference, is to growth as much as it is allowed, that is, as a line with slope $\pm1$.

More in general we have the following result.

\begin{teo}\label{teo.1.intro} 
Let  $f\in  \elle 2$,  let   $\{u_n\}\subset W^{1,p_n}(\Omega)$ be  the sequence of minimizers of $J_n$ (see \eqref{neum1}),  and let $u$ be the unique minimizer of \eqref{minlim.intro}. Then 
\[
\un \longrightarrow u \qquad \mbox{uniformly in } \overline\Omega\,.
\]
In addition, the following representations formulas hold true:
\be\label{rf}
\begin{array}{l}
\displaystyle u(x)=\max_{y\in\partial A^+} \Big[u(y)-d(x,y)\Big]\,, \quad   \forall  x\in A^+\,
\\  \mbox{ and
 } \\  \displaystyle u(x)=\min_{y\in\partial A^-}\Big[u(y)+d(x,y)\Big]\,,\quad  \forall  x\in A^-\,,
\end{array}
\ee
where  
\be
\label{Adef.intro}
A^-= \text{ supp } (u-f)^- \, ,  \qquad A^+= \text{ supp } (u-f)^+ \, ,
\ee
and $d(x,y)$ denotes the geodesic distance inside $\Omega$ (see Definition \eqref{gamma} below).  
\end{teo}

After considering the problem from a variational viewpoint, we want to address it from a PDE  perspective. Since it is not clear a priori which equetion does the minimizer of \eqref{minlim.intro} solve,   we go back to the approximating sequence of minimizers of $J_n$ (see \eqref{neum1}). The  minimizers of \eqref{neum1} are  solutions to 
\begin{equation} \label{aux.intro}
\bc
\dys u_n (x)  - \Delta_{p_n} u_n (x)  = f (x) 
\qquad & \mbox{in } \ \Omega\,,\\ 
\displaystyle 
\frac{\partial u_n}{\partial \nu } (x) =0 & \mbox{on } \ \partial\Omega \,, 
\ec
\end{equation}
where $\Delta_{p} v =$ div$(\D v |\D v |^{p-2})$, with  $p>1$,  is the   $p$-Laplacian    
and we denote  by $\nu$ the external unit normal to $\partial \Omega$. As usual, we consider such solutions in the weak sense (see Definition \ref{weakn}).
Our aim  is,  now,  to pass to the limit in \eqref{aux.intro} to find the PDE solved by the limit function, $u$.

Usually, the limit as $p \to \infty$ of solutions to equations related to the $p$-Laplacian yields to equations that involve  the   infinity Laplacian, namely,
\[
\Delta_\infty u = \langle D^2 u \nabla u, \nabla u \rangle= \sum_{i,j=1}^N u_{x_i \, x_j }\ u_{x_i  } \ u_{ x_j }\ .
\]
In particular we refer to the pioneer paper by  T. Bhattacharya, E. Di Benedetto and  J. Manfredi (see \cite{BBM})  that first studied such type of problem. We recall that  the infinity Laplacian is a second order
differential operator (in nondivergence form) that appears in many contexts. For instance,
  infinity harmonic functions (solutions to $-\Delta_\infty u
=0$) appear naturally as limits of $p$-harmonic functions
(solutions to $-\Delta_p u = -\mbox{div} (|D u|^{p-2} D
u) =0$) and they have applications to optimal transport problems, image
processing, etc. See, among the others,  \cite{acj},
\cite{BBM}, \cite{EG} and references therein.
Moreover, the infinity Laplacian plays a fundamental role in the calculus
of variations of $L^\infty$ functionals, see e.g.
\cite{ar3}, \cite{bej}, \cite{bjw1},   \cite{BLR}, \cite{cdp}, \cite{cra}, \cite{ceg},
\cite{J} and the survey \cite{acj}.
Notice that such operator is 
degenerate elliptic (that is
non-degenerate only in the direction of the gradient).

Before stating our next result,  let us introduce some notation: here and in the rest of the paper for $u, f \in C(\overline\Omega)$ we denote  
$$
\Omega^+_u=\{x\in\Omega \ : \ u(x)-f(x)>0\} \quad \mbox{and} \quad 
\Omega^-_u=\{x\in\Omega \ : \ u(x)-f(x)<0\}\,
$$
and analogously,
$$(\partial\Omega)^+_u=\{x\in\partial\Omega \ : \ u(x)-f(x)>0\}\quad \mbox{and} \quad
(\partial\Omega)^-_u=\{x\in\partial\Omega \ : \ u(x)-f(x)<0\}.$$

Now we are ready to write the problem solved by $u$. 

\begin{teo}\label{16:15}
Assume $f\in C(\overline{\Omega})$, then the unique minimizer $u \in \mathbb{X}(\Omega)$   of \eqref{minlim.intro}, is a viscosity solution to
\be\label{eqlim.intro}
\begin{cases}
 1-|\nabla u| =0 &\mbox{in} \ \Omega^+_u \,, \\
|\nabla u|-1 =0 &\mbox{in} \ \Omega^-_u \,. 
\end{cases}
\ee
Moreover it satisfies   the following boundary conditions: 
\be\label{boundarycon.intro}
\begin{cases}
1 - |\nabla  u |   \le 0 &\mbox{on} \ (\partial\Omega)^+_u \, ,\\
\displaystyle \max\left\{1 - |\nabla  u |,\frac{\partial u }{\partial \nu }\right\}\ge0 & \mbox{on } (\partial\Omega)^+_u \, ,
\end{cases}
\quad
\mbox{and}
\quad
\begin{cases}
 |\nabla  u | -1  \ge 0 &\mbox{on} \ (\partial\Omega)^-_u \, ,\\
\displaystyle \min\left\{|\nabla  u | -1,\frac{\partial u }{\partial \nu }\right\}\le0 & \mbox{on } (\partial\Omega)^-_u\,.
\end{cases}
\ee
\end{teo}
Quite surprisingly, in the  limit the second order operator is somehow lost: indeed there is no trace of the infinity Laplacian and in fact $u$ satisfies  a first order equation. 
Actually, equation \eqref{eqlim.intro} is the proper {one} in order to reflect the information carried out by \eqref{rf}: in the region $\{u-f<0\}$ the function $u$ looks like a cone pointing upwards and therefore it solves the eikonal equation (and analogously in $\{u-f>0\}$, with the eikonal equation with the reverse sign). As far as  at the boundary conditions are concerned, it is natural to wonder if it would be possible to   extend the equation in \eqref{eqlim.intro} up to $\overline{\Omega}$. However the answer is negative, and, to get convinced that \eqref{boundarycon.intro} are the natural conditions to impose on $\partial\Omega$, it is enough to look at Example 2 above.

 For the sake of clarity, notice that in Theorem \ref{16:15} $f$ is continuous and the open set $\Omega^+_u$ coincides with the interior of $A^+$.
 We will provide an explanation for the absence of second order operators in \eqref{eqlim.intro} in Proposition \ref{ifandonlyif} below.

It is also worth to mention explicitly that the regions $\Omega^+_u$ and $\Omega^-_u$ (and their boundary counterparts) depend on the solution $u$ itself. A consequence of this fact is that, even if the minimizer of \eqref{minlim.intro} is unique,   problem \eqref{eqlim.intro}-\eqref{boundarycon.intro} does not posses a   unique solution. 
In fact, we can construct minimal and maximal solutions to problem  \eqref{eqlim.intro}-\eqref{boundarycon.intro}  that,  in general, do not coincide with  the minimizer of $J_{\infty}$.

\begin{teo}\label{obstacle}
Assume that  $f \in C (\overline{\Omega})$, then   there exists a unique maximal (viscosity) solution $\overline u\in \mathbb{X}  (\Omega)$ to the following obstacle problem
\begin{equation}\label{super}
\begin{cases}
\overline u \geq  \ {f} \qquad &\mbox{in } \Omega,\\  
1 - |\nabla \overline u | \geq 0 &\mbox{in } \Omega,\\  
1 - |\nabla \overline u |   =0 &\mbox{in} \ \Omega^+_{\overline u}, \\
\end{cases}
\end{equation}
with boundary conditions
\be\label{superbc}
\begin{cases}
1 - |\nabla \overline u |   \le 0 &\mbox{on} \ (\partial\Omega)^+_{\overline u},\\
\displaystyle \max \left\{1 - |\nabla \overline u |,\frac{\partial\overline u }{\partial \nu }\right\}\ge0 & \mbox{on } (\partial\Omega)^+_{\overline u}.
\end{cases}
\ee
\end{teo}

Clearly the solution $\overline u$ provided by Theorem \ref{obstacle} solves \eqref{eqlim.intro}-\eqref{boundarycon.intro}, but it is not a minimizer of \eqref{eqlim.intro} (unless $f\in\mathbb{X}(\Omega)$). Indeed, as the intuition suggests, since the minimizer $u$ has to be close to $f$ in $L^2$, then $u-f$ has to change sign in $\Omega$ (see Figure 4 and Proposition \ref{apiuameno} for more details).

\centerline{
\includegraphics[width=0.5\textwidth]{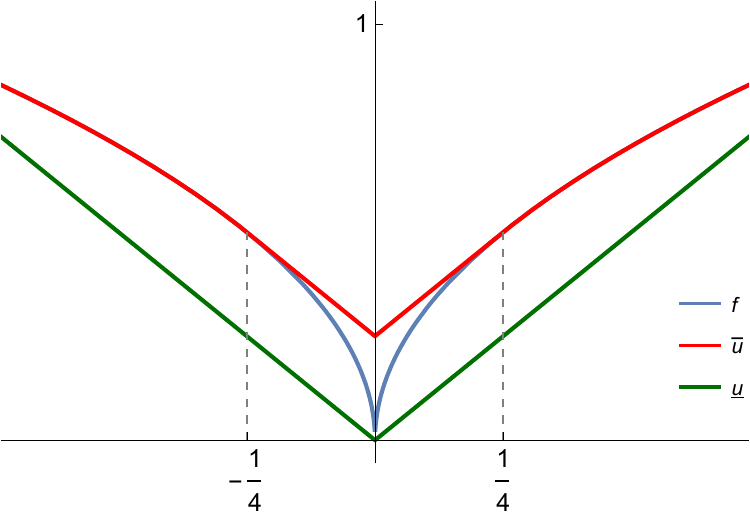}}
 \centerline{Figure 4.}

The limit procedure showed above can be, in fact, generalized to a larger class of differential equations that involve lower order terms that have a polynomial growth with respect to the gradient. More precisely, we consider the following problem 
\be \label{drift.intro}
\bc
v_n (x) - \Delta_{p_n} v_n (x) =H_n(x,\nabla v_n (x))+f(x)  
\qquad & \mbox{in } \ \Omega\,,\\
\displaystyle \frac{\partial v_n}{\partial \nu } (x) =0 & \mbox{on } \ \partial\Omega \,, 
\ec
\ee
where the Hamiltonian term  $H_n\in C(\overline{\Omega}\times\mathbb{R}^N)$ satisfies the following growth conditions
\be\label{Hgrowth}
\exists K_1,K_2>0 \,: \qquad |H_n(x,\xi)|\le p_n K_1| \xi|^{p_n}+K_2 \,, \qquad \forall x \in \Omega, \ \forall \xi \in \rn\,
\ee
and
\be\label{Hconvergence}
\exists H_{\infty}\in C(\overline{\Omega}\times\mathbb{R}^N) \, : \qquad \lim_{n\to\infty} \frac{H_n(x,\xi)}{p_n|\xi|^{p_n-4}}= H_{\infty}(x,\xi),
\ \quad \mbox{ uniformly with respect to } x \in \Omega \,. 
\ee
\begin{remark}  
Notice that the second assumption in \eqref{Hconvergence}   implies that
\be\label{Hgotozero}
\lim_{n \to \infty }H_n (x,  \xi)= 0 \quad \mbox{for any } |  \xi| <1\,, \quad \mbox{uniformly with respect to  } x \in \Omega.
\ee
\end{remark}

Among all the possible choice of $H_n$ that satisfy \eqref{Hgrowth}, let us mention two important  examples: 
\[
H_n(x,\xi)=p_n |\xi|^{p_n-2}\xi\cdot B(x)  \quad \mbox{ with }\quad B\in (C(\overline{\Omega}))^N\ \quad  \mbox{and} \ \quad  \  H_n(x,\xi)=p_n |\xi|^{p_n} \,.
\]
 These two types of lower order terms, accounting for the ``linear $(p-1)$ growth" and ``natural growth", are thoroughly studied in   literature (for  fixed $p$). Typically, the presence of the Hamiltonian affects the coercivity of the operator and some specific techniques are required to obtain apriori estimates: the first case can be dealt with symmetrization procedures \cite{bfm}, $log$-type estimates \cite{lucio}, slicing techniques \cite{abilitazione} or arguing by contradiction \cite{nsesamai}. For the case of natural growth we mention the exponential test function approach developed in \cite{BMP1}--\cite{BMP}, again symmetrization in \cite{ALT}, and an argument based on Schauder fixed point Theorem in \cite{MN}.
 
Let us notice that, under assumption \eqref{Hgrowth}, for a smooth $f$ and a fixed $n$, the existence of a solution $v_n\in W^{1,p_n} (\Omega) \cap \lio$ to problem  \eqref{drift.intro} is a straightforward consequence of the results contained in  \cite{BMP} (see Theorem \ref{exn} below for more details).  We also explicitly point out that the presence of the zeroth lower order term allows us to solve \eqref{Hconvergence} without any smallness assumption on the datum $f$.

The mail difficulty when one wants to pass to the limit as $p$ diverges  in \eqref{drift.intro}  is to obtain apriori estimates for $v_n$ that are stable with respect to $n$. Up to our knowledge, the only related result for quasilinear equations with gradient lower order terms is contained in \cite{BL}, where a large solution problem is considered. 
 
Our main result concerning problem \eqref{drift.intro} is the following. 
 
\begin{teo} \label{33}
Assume \eqref{Hgrowth}--\eqref{Hconvergence} and let $f\in C(\overline{\Omega})$. Then there exists $v\in W^{1,\infty}(\Omega)\cap \lio$ such that
\[
v_n \to v\quad \mbox{uniformly in } \overline {\Omega}
\] 
which solves (in viscosity sense)   the following   equation
\be\label{driftinftyIntro}
\begin{cases}
\max\big\{1-|\nabla v|,-\Delta_{\infty}v- H_{\infty}(x, \nabla v)   \big\}=0 \qquad &\mbox{in} \quad  \Omega^+_v\,,\\
\min\big\{|\nabla v|-1,-\Delta_{\infty}v- H_{\infty}(x, \nabla v) \big\}=0 \qquad &\mbox{in} \quad  \Omega^-_v\,\,.
\end{cases}
\ee
together with the following boundary conditions:
\be\label{bc1}
\begin{split}
\min\left\{\max\big\{1-|\nabla v |,-\Delta_{\infty}v- H_{\infty}\big(x,\nabla v \big)     \big\} , \frac{\partial  v}{\partial \nu}\right\}\le 0 \quad \mbox{on } \   (\partial\Omega)^+_v \,,\\
\max\left\{\max\big\{1-|\nabla v|,-\Delta_{\infty} v - H_{\infty} \big(x,\nabla v \big)     \big\} , \frac{\partial v}{\partial \nu}\right\}\ge 0 \quad \mbox{on } \    (\partial\Omega)^+_v \,, \\
\end{split}
\ee
and
\be\label{bc2}
\begin{split}
 \min\left\{\min\big\{|\nabla v|-1,-\Delta_{\infty} v- H_{\infty}\big(x,\nabla v\big)     \big\}, \frac{\partial v}{\partial \nu}\right\}\le 0 \quad \mbox{on } \     (\partial\Omega)^-_v\, ,\\
  \max\left\{\min\big\{|\nabla v |-1,-\Delta_{\infty}v - H_{\infty} \big(x,\nabla v\big)     \big\} , \frac{\partial v}{\partial \nu}\right\}\ge 0 \quad\mbox{on } \    (\partial\Omega)^-_v\,. 
\end{split}
\ee
\end{teo}

At the first glance we notice that the structure of \eqref{driftinftyIntro} is more involved then the one in \eqref{eqlim.intro}. First of all $v$ does not belong in general to $\mathbb{X} (\Omega)$ since the bound on $\nabla v$ depends also on the Hamiltonian (see Section \ref{lowerorderterm}). 
Moreover, as expected, the presence of the Hamiltonian term changes both  the limit equation and the boundary conditions.

Let us now go back to the minimizer $u\in \mathbb{X} (\Omega)$ of \eqref{minlim.intro} and notice that it also solves the equation \eqref{driftinftyIntro} with $H_{\infty}\equiv0$. Apparently, this provides some more information about $-\Delta_{\infty}u$ that is missing in \eqref{eqlim.intro}. However, we have the following result.

\begin{prop}\label{ifandonlyif}
Assume that  $f\in C(\overline{\Omega})$,  and that $w\in \L1$. Then $w$ solves \eqref{eqlim.intro} if and only if it solves 
$$
\begin{cases}
\max\big\{1-|\nabla w|,-\Delta_{\infty}w \big\}=0 \qquad &\mbox{in} \quad  \Omega^+_w \,,\\
\min\big\{|\nabla w|-1,-\Delta_{\infty} w  \big\}=0 \qquad &\mbox{in} \quad  \Omega^-_w\,\,.
\end{cases}
$$
\end{prop}

Roughly speaking this is due to the bound $|\nabla w|\le1$ and the fact that cones are  solutions to   the eikonal equation and are also infinity harmonic away from their singular points.

Finally, we include some considerations about a very well known family of functionals that allow   \emph{free discontinuities}. Unfortunately, we are able to provide only partial results in this case and we left a more complete analysis of the limiting behavior of such functionals for $p\to\infty$ for future research. Let us consider the functional
$ \mathcal{I}_{p}:SBV(\Omega)\to \mathbb{R}\cup\{+\infty\}$ defined as
\be\label{ms1}
 \mathcal{I}_{p}(v)=\frac1{p}\io|\nabla v|^{p} +\frac1{r}\io|v-f|^r+\mathcal{H}^{N-1}(S_v)\,,
\ee
with $r\geq 1$. 
Here $SBV(\Omega)$ is the space of Special Functions of Bounded Variation, $\nabla v$ is the absolutely continuous part of $Dv$ with respect to the Lebesgue measure, $S_v$ is the singular set of $v$ and $\mathcal{H}^{N-1}$ is the $N-1$ dimensional Hausdorff measure (more details and formal definitions can be found in Section \ref{mumford}). We 
set $ \mathcal{I}_{p}(v) =+\infty$ when $\nabla v \not\in L^{p} (\Omega)$.

The functional  $ \mathcal{I}_{p }(v)$  was introduced in the seminal paper by  De Giorgi, Carriero and Leaci (\cite{DCL}) to provide a weak formulation of the Mumford–Shah image segmentation problem, namely, minimizing
\be\label{mumsha}
F_{p} (K,v)=\frac1{p}\io|\nabla v|^{p}+ \frac1r\io|v-f|^r+\mathcal{H}^{N-1}(K),
\ee
among all $K\subset\overline\Omega$ closed sets,  with  $v\in W^{1,p}(\Omega\setminus K)$.
An in-depth overview about the history of these functionals, their relation with applications, and the huge impact that their study has had in the field of Calculus Variation is definitively out of reach for our contribution. We refer the interested reader to the original paper \cite{DCL}, the classical monograph \cite{AFP}, and to the nice review \cite{fusco}. For our aims, we simply recall that the main difficulty in dealing directly with \eqref{mumsha} is that the Hausdorff measure is not lower semicontinuous with respect to any reasonable metric the collection of closed set of $\Omega$ can be equipped with. Therefore, the strategy of \cite{DCL} was to find a minimizer $u\in SBV(\Omega)$ of $ \mathcal{I}_{p}$ and then show that $(\overline{S_{u}},u)$ was also a minimizer to \eqref{mumsha}. The most delicate step in their argument was to prove that $\mathcal{H}^{N-1}(\overline{S_{u}}\cap\Omega)=\mathcal{H}^{N-1}({S_{u}})$ and this crucial step was achieved by a lower bound on the density of $S_{u}$, namely
\be\label{lowerdensity}
\exists \theta >0 \: \quad \mathcal{H}^{N-1}(S_{u}\cap B_{\rho}(x))>\theta \rho^{N-1}, \qquad \forall x\in \overline{S_{u}}, \, \forall \rho>0\,
\ee
(see Theorems 7.15, 7.21 and 7.22 in \cite{AFP}). 

At this point it is natural to wonder what happens with the limiting minimization problem relative to \eqref{mumsha} as $n\to\infty$. The answer is given by the following Proposition.

\begin{prop} 
Assume that $f\in \lio$ and that $p_n\to +\infty$, then 
the sequence of functionals $ \mathcal{I}_{p_n} :SBV(\Omega)\to \mathbb{R}\cup\{+\infty\}$ defined in \eqref{mumsha} $\Gamma$-converges to $ \mathcal{I}_{\infty}:SBV(\Omega)\to \mathbb{R}\cup\{+\infty\}$ defined as
\begin{equation}\label{Minfty}
 \mathcal{I}_{\infty}(v)=
\begin{cases}\displaystyle
\frac1r\io|v-f|^r+\mathcal{H}^{N-1}(S_v) & \mbox{ if } \ |\nabla v|\le1\,,\\
+\infty & \mbox{otherwise}.
\end{cases}
\end{equation}
\end{prop}

The next step would be to prove that a minimizer of \eqref{Minfty} is also a minimizer to 
\be\label{mumshainf}
F_{\infty}(K,u)=\frac1r\io|u-f|^r+\mathcal{H}^{N-1}(K),
\ee
with $K\subset\overline\Omega$ a closed set, $u\in W^{1,\infty}(\Omega\setminus K)$ and $|\nabla u|\le 1$. Unfortunately, we are not able to show such a property since we do not know how to obtain a lower density bound, as in \eqref{lowerdensity}, for the limit function $u$. We leave this as an interesting open problem.

\section*{Notations and basic definitions}

Let us introduce the notion of geodesic distance (see \cite{rs}, for instance):  for $x,y\in\Omega$  we define 
\begin{equation}\label{gamma}
\begin{array}{l}
\displaystyle
\dys \Gamma= \{  \gamma_{x,y} \in C^1 ( [0,1] ; \Omega)  \,  \mbox{ with } \,  \gamma(0)=x, \ \gamma(1)=y   \} \\[8pt]
\quad 
\mbox{ and}  
\quad \\[8pt]
\displaystyle
d(x,y)=\inf_{\gamma \in \Gamma} \left\{\mathcal{L} (\gamma_{x,y} )   \right\}, \quad \mbox{ where  } \ 
\mathcal{L} (\gamma_{x,y} ) =  \mbox{ lenght } (\gamma_{x,y} )\,.
\end{array}
\end{equation}
Of course, if $\Omega$ is convex $d(x,y)\equiv|x-y|$, but for general domains we just have that $d(x,y)\ge |x-y|$. The geodesic distance is the natural quantity to consider in our setting since the condition 
$$v\in \mathbb{X} (\Omega) := \{ v \in W^{1,\infty} (\Omega) \,:\, |Dv | \leq 1 \ \ \ a.e. \ x \in \Omega\}
$$ can be rephrased as 
\[
v \in \left\{ v \in W^{1,\infty} (\Omega) \,:\,\, |v(x)-v(y)|\le d(x,y), \quad \forall \ x,y\in \Omega \right\} .
\]

Now, let us introduce some extra notations that will be used in the rest of the paper.
Let us recall that (see Proposition 4.17 in \cite{Br}) the support of a general function $f:\Omega\to\mathbb{R}$  is defined as
\[
\mbox{supp}(f)=\left(\bigcup_{i\in I}\omega_i\right)^c
\qquad \mbox{ 
where
 }
\qquad I=\Big\{\omega \mbox{ open subset of } \ \Omega \ \mbox{such that}\ f=0 \ a.e. \ \mbox{ in }\omega\Big\}.
\]
By definition $\mbox{supp}(f)$ is a closed set with respect to the subset topology relative to $\Omega$. For any $f\in\elle2  $, and for any continuous function $u$ we define the following sets
\be
\label{Adef}
A^-= \text{ supp } (u-f)^- \, ,  \qquad A^+= \text{ supp } (u-f)^+ \,.
\ee

For future use we need to give the precise definition of a solution to a partial differential equation in the viscosity sense.

\begin{defin} \label{visc} Let $D$ be a locally compact subset of $\mathbb{R}^{\mathbb{N}} $ and let  $F:D\times\mathbb{R}\times\mathbb{R}^{\mathbb{N}} \times \mathcal{S} (N)$ be a continuous function that satisfies
\[
F(x,r,\xi,X)\le F(x,s,\xi,Y) \ \ \ \mbox{whenever} \ \ \ r\le s \ \ \ \mbox{and} \ \ \ Y\le X .
\]
We say that $u\in C(D)$ is a viscosity sub (super) solution  to
\be\label{f}
F(x,u,\nabla u, D^2 u)=0 \qquad \mbox{  at } \ x \in D
\ee
if, for any $C^2(D)$ function $\varphi$ such that $u-\varphi$ has a strict local maximum (minimum) at $x_0$,  it holds true
\[
F(x_0,\varphi(x_0),\nabla \varphi(x_0), D^2  \varphi(x_0) )\le (\ge) \ 0.
\]
Finally, $u\in C(D)$ is a viscosity solution of \eqref{f} if it is both sub and super solution.
\end{defin}

Let us denote by $\delta (x) = $ dist$(x,\partial \Omega)$ and let us recall that $\delta (x) \in W^{1,\infty}_0 (\Omega)$ and that it satisfies (in the viscosity sense) the eikonal equation 
$$
| \nabla  \delta |-1=0 \quad \mbox{ in } \Omega\,.
$$

In the sequel we will use the following result.

\begin{prop}\label{L1}
Any function in $\L1$ satisfies  both 
\be \label{doublel}
|\nabla u|-1 \le 0 \qquad \mbox{and} \ \ \ -|\nabla u|+1\ge0 \qquad \mbox{in} \  \Omega
\ee
in the viscosity sense. 
\end{prop}

\begin{proof}[Proof of Proposition \ref{L1}]
Let us observe that   \eqref{doublel} is equivalent to say that the following inequalities hold true at any $x_0\in \Omega$:
\begin{itemize}
\item for any $\vp \in C^2(\Omega)$ such that there exists $r>0$ with $u(x) - \vp (x)  \leq u(x_0) - \vp (x_0) $ in $B_{r} (x_0) \subset \Omega$, we  have that $
|\nabla \vp (x_0) |-1 \le 0$; 

\item for any $\vp \in C^2(\Omega)$ such that there exists $r>0$ with $u(x) - \vp (x)  \geq u(x_0) - \vp (x_0) $ in $B_{r} (x_0)\subset \Omega$,   we have that $
-|\nabla \vp (x_0) |+1 \ge 0$. 
\end{itemize}

 In order to prove the first inequality in \eqref{doublel} let $x_0\in\Omega$ and assume that there exists $\varphi$, a $C^2$ function of a neighborhood of $x_0$, such that $u-\varphi$ has a strict local maximum at $x_0$, with  $u(x_0)=\varphi(x_0)$. 
 Suppose by contradiction that $|\nabla \varphi(x_0)|>1$. Taking the first order Taylor expansion of $\varphi(x)$ at $x_0$  it follows that
\[
\varphi(x)-\varphi(x_0)=  \nabla \varphi (x_0)(x-x_0)+o(|x-x_0|)  \ \ \ \mbox{as} \ \ \ x\to x_0.
\] 
Let   us consider 
$$x_t=- \frac{\nabla \varphi (x_0)}{|\nabla \varphi (x_0)|} t+x_0, \quad \mbox{ for }  t>0
\mbox{ small}.$$
We have  $|x_t - x_0|= t $ and thus 
\[
u(x_t)-u(x_0)\le \varphi (x_t)-\varphi(x_0) 
=   \nabla \varphi(x_0) \cdot (x_t-x_0) +o(|x-x_0|) = 
 - |\nabla \varphi(x_0)| t \big(1+o(1)\big)\, ,  \ \ \ \mbox{as} \ t\to0\,.
\]
This is contradicts   $u\in \mathbb{X} (\Omega) $. 

The second inequality in \eqref{doublel} follows in the same way. 
\end{proof}

\section{The best Lipschitz approximation with gradient constraint} 

Let $p_n$ be an increasing sequence  such that 
$p_n\to \infty$ as $n\to\infty$.  Let us set, for any $f\in \elle2$,  the following functional 
\be\label{neum}
J_{n}(v)=\frac1{p_n}\io|\nabla v|^{p_n}+\frac12\io(v-f)^2 \,, \ \ \qquad \ v\in W^{1,p_n} (\Omega)\,.
\ee
For any fixed $n$, it is hot hard to see that there exists a unique (exploiting the convexity of \eqref{neum}) minimizer  $u_n\in W^{1,p_n} (\Omega)$  to \eqref{neum}. 
Thus, from direct variational arguments, the minimizer, $u_n$, turns out to be a weak solution to    
\begin{equation} \label{aux}
\bc
u_n (x) - \Delta_{p_n} u_n (x) = f(x)  
\qquad & \mbox{in } \ \Omega\,,\\
\displaystyle \frac{\partial u_n}{\partial \nu } (x) =0 & \mbox{on } \ \partial\Omega \,. 
\ec
\end{equation}
Let us recall here such a formulation. 

\begin{defin}\label{weakn}
A weak solution    $u_n$ to \eqref{aux} is a $ W^{1,p_n} (\Omega)$ function that  satisfies 
\begin{equation} \label{auxw}
\io u_n (x) \varphi + \io |\nabla u_n (x)|^{p_n-2} \nabla u_n \cdot \nabla \varphi = \io f \varphi   \,,
\qquad  
\forall \varphi \in W^{1,p_n} (\Omega)\,.
\end{equation}
\end{defin}

The aim of this section is to deal with the limiting behaviour of these three objects - the minimizer, the functional, and the equation - as $p_n$ diverges.

 We start with the following proposition that provides  a limit for the sequence of minimizers together with a limiting minimization problem associated to the limit of \eqref{neum}.

\begin{prop}\label{auxlem} Assume $f\in\elle 2$ and let $\{u_n\}\subset W^{1,p_n}(\Omega)$ be the sequence of minimizers of \eqref{neum}, then there exists $u \in \mathbb{X}(\Omega)$ such that
\[
\un \longrightarrow u \qquad \mbox{uniformly in } \overline\Omega .
\]
Moreover, $u$ is the unique minimizer of  the functional 
\begin{equation}\label{16:50}
J_{\infty}(v)=\frac12\io(v-f)^2 \,,\ \ \ v\in \mathbb{X}(\Omega).
\end{equation}
\end{prop}

\begin{proof}
The minimality of $\un$ gives us that
\[
\frac{1}{2}\io(u_n-f)^2+\frac1{p_n}\io|\nabla \un |^{p_n}\le J_{n}(v)\le J_{n}(0)=\frac12\io|f|^2.
\]
This implies that
\be\label{estimates}
\|u_n\|_{\elle 2}\le \|f\|_{\elle 2} \qquad \mbox{and} \qquad \|\nabla u_n\|_{\elle{p_n}}\le   \big(p_n \|f\|_{\elle 2}^{2} \big)^{\frac1{p_n}}.
\ee
For $N < p_n$, let us recall that Morrey inequality (see \cite[Corollary 9.14]{Br}) implies that, for any $x,y \in \Omega$,
\be\label{morrey1}
|u_n(x)-u_n(y)|\leq C_M \|\nabla u_n\|_{\elle {p }} ,
%\le |v(y)|+4\mbox{diam}(\Omega)\|\nabla v\|_{\elle {p_n}} , \quad \forall \ x,y\in \Omega
\ee
where $C_M=C_M(N,\Omega)$ can be chosen independent on $p$ (see in particular formula (28) in the proof of  \cite[Theorem 9.12]{Br}).
Therefore, we deduce that
\be\label{morrey}
|u_n(x)|\le |u_n(y)|+{C}_M \|\nabla u_n\|_{\elle {p_n}} , \quad \forall \ x,y\in \Omega\,,
\ee
and consequently
\be\label{morreybis}
\|u_n\|_{\elle{\infty}}\le |\Omega|^{\frac12}\|u_n\|_{\elle 2}+\tilde{C}_M \|\nabla u_n\|_{\elle {p_n}}.
\ee
Combining 
 \eqref{estimates}  with \eqref{morreybis} we conclude that
\be\label{limitate}
\|u_n\|_{\elle{\infty}}\le |\Omega|^{\frac12}\|f\|_{\elle 2}+c_\Omega \|f\|_{\elle 2}^{\frac2{p_n}}\,.
\ee
On the other hand, using H\"older's inequality, it follows that, for any $q$ such that $N<q<p_n$,
\[
\|\nabla u_n\|_{\elle q}\le \|\nabla u_n\|_{\elle {p_n}}|\Omega|^{\frac1q-\frac 1{p_n}}\le  p_n^{\frac1{p_n}}
 \|f\|_{\elle 2}^{\frac2{p_n}}|\Omega|^{\frac1q-\frac 1{p_n}}.
\]

Thanks to the two estimates above, we deduce that the sequence $\{\un\}$ is bounded in $W^{1,q}(\Omega)$ and then up to a (not relabeled) subsequence, $\un\to u$ weakly in $W^{1,q}(\Omega)$ and uniformly in $\Omega$. The lower semi continuity of the norm also implies that for any $q\geq 1$ we have
\[
\|\nabla u \|_{\elle{q}}\le |\Omega|^{\frac1q}.
\]
Passing to the limit as $q\to \infty$ we conclude that $u\in \mathbb{X}(\Omega)$.

In order  to show that $u$ is, in fact, a minimizer to \eqref{16:50}, let us notice that, by definition, $u_n$ satisfies 
\[
\frac12\io(u_n-f)^2\le \frac12\io(u_n-f)^2+  \frac{1}{p_n}\io|\nabla u_n|^{p_n}\le  \frac12\io(z-f)^2+\frac{1}{p_n}\io|\nabla z|^{p_n}\,,  
\]
$ \forall z\in\mathbb{X}(\Omega)\subset W^{1,p_n}(\Omega)$.
Since $|\nabla z| \leq 1$ a.e. in $\Omega$ we get
\[
\frac12\io(u_n-f)^2
\le  \frac12\io(z-f)^2+\frac{1}{p_n} | \Omega | \,, \quad \forall z\in\mathbb{X}(\Omega)\subset W^{1,p_n}(\Omega)\,.
\]
We can now pass to the limit in the previous inequality (recall that $u_n \to u$ uniformly in $\Omega$), and we deduce that $u$ is a minimizer. 
The uniqueness follows by a  the convexity of $J_\infty$.
\end{proof}

\begin{remark}  
When $f \in L^\infty (\Omega)$ we have a explicit uniform bound for $u_n$, $$|\un|\le \|f\|_{\elle{\infty}}=\vcentcolon \overline k .$$ Indeed, $T_{\overline k}(\un)$ is a competitor for $\un$ for the energy \eqref{neum} and the uniqueness of the minimizer of $J_n$ implies $\un \equiv  T_{\overline k}(\un)$. 
\end{remark}

In other words, the minimizer of \eqref{16:50} is the best  approximation of $f$ among all the   functions in $\L1$. The gradient constraint imposes to the minimizer a specific geometric structure. The goal of the following results is to clarify such a  structure.

\medskip

\begin{prop}\label{apiuameno}
Assume that $f\in\elle2$ and  let $u\in\mathbb{X}(\Omega)$ be the minimizer to \eqref{16:50}.  Then 
$$|A^-||A^+|=0 \qquad \mbox{  if and only if } \qquad f\in \mathbb{X} (\Omega)\,.$$
\end{prop}
\begin{proof} We observe that if $f\in \L1$, then $u \equiv f$
and hence $|A^-|=0$ and $|A^+|=0$. Suppose now, that $f\in \elle2 \setminus \L1$ and 
by contradiction suppose that $|A^-||A^+|\neq 0$. If both $|A^-|$ and $|A^+|$ are zero then
$u=f$ and we reach a contradiction. 

Assume now  that $|A^-|>0$ and $|A^+|=0$
(the case $|A^+|>0$ and $|A^-|=0$ is analogous). 
Then, for any $\eps>0$,  $u+\eps \in \L1$  satisfies 
\[
 J_{\infty}(u+\eps)= J_{\infty}(u)+\frac12 \eps^2|\Omega|+ \eps\int_{A^-}(u-f)= J_{\infty}(u)+\frac12 \eps\left(\eps|\Omega|-2\int_{A^-}|u-f|\right),
\]
and the last term in the right hand side is negative for $\eps$ small enough. This contradicts the minimality of $u$. 
\end{proof}

\begin{remark}
{\rm In fact, the proof shows that 
$$|A^-|=0 \mbox{ and } |A^+|=0 \qquad \mbox{  if and only if } \qquad f\in \mathbb{X} (\Omega)\,.$$
However, note that the statement of the proposition highlights that 
it is enough to have $|A^-|=0$ or $ |A^+|=0$ to obtain that $u=f\in \mathbb{X} (\Omega)$.}
\end{remark}

Next we deduce a representation formula for the unique minimizer to \eqref{16:50}.

\begin{teo}\label{minmaxcone} 
Assume $f\in  \elle 2$ 
and let  $u$ be the unique minimizer to \eqref{16:50}, then  the following representations formulas hold true:
\be\label{rf1}
\begin{array}{l}
\displaystyle u(x)=\max_{y\in\partial A^+} \Big[u(y)-d(x,y)\Big]\,, \ \ \ \forall  \ x\in A^+\,
\end{array}
\ee
 and
 \be\label{rf2}
\begin{array}{l}
\displaystyle u(x)=\min_{y\in\partial A^-   }\Big[u(y)+d(x,y)\Big]\,,\ \ \ \forall  \ x\in A^-.
\end{array}
\ee
\end{teo}

\begin{proof}
We assume that $f\notin \mathbb{X}    (\Omega) $, otherwise $u\equiv f$ and there is nothing to prove. Proposition \ref{apiuameno} assures us that both $A^+$ and $A^-$ are nontrivial, then the max and the min in \eqref{rf1} and \eqref{rf2} are well defined.
We start by proving the second formula above. Let us define 
$$
v(x)=
\begin{cases}
\dys \min_{y\in\partial A^- }\Big[ u(y)+d(x,y)\Big]\,\qquad &\mbox{ in } \quad A^-\,,
\\[1.5 ex]
  u(x)   &\mbox{ in } \quad \Omega \setminus A^- \,.
\end{cases} 
$$
The strategy of the proof is to show that such a $v$ belongs to  $\mathbb{X}  (\Omega)$ and that it is a competitor for $u$ in the minimization of \eqref{16:50}. Then, uniqueness
of the minimizer   implies that
$$v\equiv u.$$ 

The first step is to show that 
\be\label{partiallip}
|v(x)-v(z)|\le d(x,z)\,, \qquad \forall \ x,z\in A^-.
\ee
By definition we have that
\be\label{lip1}
\begin{array}{l}
\displaystyle 
v(x)=u(\bar x)+d (x,\bar x)=\min_{y\in\partial A^-   } \Big[u(y)+d(x,y)\Big] 
\\[8pt]
\qquad\mbox{and} \qquad \\[8pt]
\displaystyle v(z)=u(\bar z)+d (z,\bar z)=\min_{y\in\partial A^-  }\Big[ u(y)+d(x,y)\Big]\,,
\end{array}
\ee
for some  $\bar x, \bar z\in  \partial A^- $. If $\bar x=\bar z$ we have that
\[
|v(x)-v(z)|=d(x,\bar x)-d(z,\bar x)\le d(x,z).
\]
On the other hand, if $\bar x \neq \bar z$, we get
\be\label{lip2}
\begin{array}{l}
\displaystyle
v(x) = u(\bar x)+d(x,\bar x)\le u(\bar z)+d(x,\bar z) \\[8pt]
\qquad \mbox{and} \qquad \\[8pt]
\displaystyle 
v(z)= u(\bar z)+d(z,\bar z)\le u(\bar x)+d(z,\bar x).
\end{array}
\ee
Thus by  \eqref{lip1}--\eqref{lip2}, we deduce that  
\[
-d(x,z)\le d(x,\bar x)-d(z,\bar x)\le v(x)-v(z)\le d(x,\bar z)-d(z,\bar z)\le d(x,z),
\]
and \eqref{partiallip} follows. 

Next we prove that
$$
v\equiv u \quad \mbox{on } \partial A^- . $$ 
By contradiction, assume that there exists  $x\in \partial A^- 
 $ such that, say,  $u(x) >v(x)$. Thus,  $\exists \bar x\in \partial A^-\,, \, x \neq \bar x$ such that   
$ v(x) =u(\bar x)+d(x,\bar x)< u(x) $ that 
  contradicts  $u\in \mathbb{X}  (\Omega)$.

We need to prove now that $v\in \mathbb{X} (\Omega)$. The only relevant case we have to consider is when $y\in \Omega \setminus A^-$ and $x\in A^- $.  Consider   $\gamma \in \Gamma$  (see \eqref{gamma} for the definition),   
let $t^* := \inf \{t \in (0,1) \,: \, \gamma(t) \not\in \Omega\}$, and $z= \gamma (t^*)$.   
Thus, since $z \in \partial A^-$, then $u(z)=v(z)$ and 
 \[
|v(x)-v(y)|\le |u(x)-u(z)|+|v(z)-v(y)|\le d(x,z)+d(z,y)\le \mathcal{L} (\gamma_1)+ \mathcal{L}(\gamma_2)= \mathcal{L} (\gamma), 
\]
where $\gamma_1 = \gamma \mres\big( t \in [0,t^*]\big)$ and 
$\gamma_2 = \gamma \mres  \big( t \in [t^*,1]\big)$. 
Minimizing the above expression with respect to any $\gamma\in \Gamma$ we deduce that  $v\in \mathbb{X} (\Omega)$.

To conclude the proof we need to show that $v\equiv u$ in $A^-$. If we assume that there exists $\bar x\in A^- $ such that $u(\bar x)> v(\bar x)$, then, by definition of $v$, we have that $v(\bar x)=u(\bar y)+d(\bar x, \bar y)$ for some $\bar y\in\partial A^-$. It follows that 
\[
u(\bar x)-u(\bar y)>v(\bar x)-v(\bar y)=d(\bar x,\bar y),
\]
that contradicts the fact that $u\in \L1$.
Therefore, let us assume, arguing again by contradiction, that there exists $\bar x\in A^-$ such that 
$$\dys \sup_{y\in A^-} (v-u)(y) = (v-u)(\bar x)=\mu >0\,. $$
We set 
\[
D=\bigg\{x\in A^- \ : \ v(x)-u(x)>\frac{\mu}2\bigg\}\,,
\quad \dys 
\eta= \essinf_{y\in D} \ (f-u)  \geq 0\  
\]
 and we define
\[
\tilde u= \begin{cases}
\max\big\{u,v-\varepsilon \big\} &\mbox{in} \ D\, ,\\
u & \mbox{in} \ \Omega\setminus D \,, 
\end{cases}
\]
 with  
$\varepsilon=\max\left\{\mu-\frac{\eta}{2}, \frac{\mu}2 \right\}$. 
 Since both $u$ and $v$ are continuous 
it follows that $\tilde u$ is continuous in $\Omega$, and moreover $\tilde u\in \L1$ (since both $u$ and $v$ are).  
Observe now that that $\tilde u$ satisfies
 $\tilde u <f$ in $A^-$: indeed in $A^- \setminus D$ this follows by the definition of $A^-$, while
 $$ \tilde u (x) -f(x) = v(x) - \eps -f (x)= \underbrace{v(x) - u(x)}_{\leq{\mu} } - \eps + \underbrace{u(x) - f (x) }_{\leq - \eta } \leq \mu - \eps - \eta <0\qquad \mbox{in } D\,.$$
 Consequently 
\be\label{competitor}
\big(u(x)-f(x)\big)^2\ge \big(\tilde u (x)-f(x)\big)^2 \ \ \ \mbox{in} \ \Omega, 
\ee
since the two quantities above coincide in $\Omega\setminus D$  while  
\[
u(x)-f (x)< \max\{u(x),v(x)-\varepsilon\}-f(x)< 0 \qquad \mbox{ in } D.
\]
Then, we have that  
$$J_{\infty} (\tilde u ) \ <\  J_{\infty} (u) $$ 
that is in contradiction with the fact that $u$ minimizes $J_{\infty}$. 

Thus, we have that $v(x) = u(x) $ in $\Omega$ and in particular  $v(x) = u(x) $ in $A^-$ and consequently  
  in $A^-$ the representation formula given by  \eqref{rf2} holds true. 

Analogously one can prove that also   \eqref{rf1} is in force. 
\end{proof}

\begin{proof}[Proof of Theorem \ref{teo.1.intro}]
The proof follows immediately putting together the results of Proposition \ref{auxlem} and Theorem \ref{minmaxcone}.
\end{proof}

% \begin{remark}
%We finally observe that the procedure that we have studied in the previous sections can be performed even for functionals defined in $W^{1,p_n}_0 (\Omega)$. The only difference that we can find is that the boundary condition now are the ones that follows from such a constraint, i.e. 
%
% \end{remark}

\begin{remark}
Let us observe that a similar version of the previous results can be obtained working in $W^{1,p_n}_0 (\Omega)$, i.e.  considering homogeneous Dirichlet boundary conditions. Actually, this latter case can be recovered as a sort of corollary of the $W^{1,p_n}(\Omega)$. To get convinced of this fact notice that the limit  space for the Dirichlet problem is
$$
\mathbb{X}_0 (\Omega) := \Big\{ v \in W^{1,\infty} (\Omega) \,:\, |Dv | \leq 1 \ \ a.e. \mbox{ in } \Omega \mbox{ and } u=0 \mbox{ on } \partial \Omega \Big\}\,.
$$
We consider the two minimization problems
\[
\inf_{\mathbb{X}_0(\Omega) }  \frac12 \io(v-f)^2 \qquad \inf_{\mathbb{X}(\Omega) }  \frac12 \io(v-\tilde{f})^2
\]
where  $\tilde{f}(x)=\max\{-\delta(x),\max\{f(x),\delta(x)\}\}$. We claim that the minimizer $u\in \mathbb{X}_0(\Omega) $ of the first functional  coincides with the minimizer $\tilde u\in \mathbb{X}_0(\Omega) $ of the second one. 
Indeed, thanks to the definition of $\tilde{f}(x)$ one easily deduce that $\tilde u\in \mathbb{X}_0(\Omega)$. Moreover if $u$ and $\tilde u$ would differ in some region one could build a competitor for one of the two, contradicting uniqueness.
 \end{remark}

\section{The problem with gradient lower order term }\label{lowerorderterm}

In this section we address the limit, as $n\to\infty$, of the family of problems   \eqref{drift.intro}.

Our first result provides existence and uniqueness of a solution for problem \eqref{drift.intro}, for any fixed $p_n\in(1,\infty)$, together with some  estimates useful in order to study the asymptotic 
behaviour as $n\to\infty$.
\begin{teo}\label{exn}
Let us assume \eqref{Hgrowth} and $f\in \elle{\infty}$. Then, for any fixed $n\in \mathbb {N}$, there exists a unique solution $v_n \in W^{1,p}(\Omega) \cap \lio $ to problem \eqref{drift.intro} in the following weak sense
\begin{equation}\label{drift}
\io v_n \phi+\io |\nabla v_n|^{p_n-2}\nabla v_n\nabla \phi=\io \big(H_n(x,\nabla v_n)+f\big)\phi \,,\quad \forall \phi\in W^{1,p_n}(\Omega)\cap \elle{\infty}.
\end{equation}
Moreover the following estimates hold true
\be\label{pestimate}
\|v_n\|_{\elle{\infty}}\le A :=  K_2+\|f\|_{\elle{\infty}} \quad \mbox{and} \quad \|\nabla v_n\|_{\elle{p_n}}\le B_{n} :=  A^{\frac{1}{p_n}} e^{ (K_1+1/p_n)A}  |\Omega |^\frac1{p_n}\,.
\ee
\end{teo}

\begin{remark}\label{limitB}{\rm
Observe that $$
\lim_{n\to +\infty} B_n = e^{K_1 (K_2+\|f\|_{\infty})}=\vcentcolon B_{\infty}>1\,.
$$
}
\end{remark}

\begin{proof} The existence of a bounded  weak solution to \eqref{drift} follows by \cite{BMP1} -- \cite{BMP}, while the uniqueness is a consequence of  \cite[Theorem 1.2]{LPR}. 

In order  to prove  that $\| v_n\|_{\elle{\infty}}\le A$, assume by contradiction that
\[
\mbox{meas}(\{x\in \Omega\,: \ \  |v_n|> A\}) >0\,;
\]
therefore, there exists $k>A$ such that $$\mbox{meas}(\{x\in \Omega\,: \  |v_n|>k\})>0.$$ 
Let us consider $\mu_n=p_n K_1+1$ an let us choose for $k\geq 0$,  
 $$\phi = (e^{\mu_n  |G_k(v_n )|}-1)\mbox{sign}(v_n )$$ as a test function in \eqref{drift}. Using assumption \eqref{Hgrowth}, we get that
\[
\begin{split}
  \io |v_n|(e^{\mu_n  |G_k(v_n)|}-1)+\mu_n\io|\nabla G_k(v_n)|^{p_n} e^{\mu_n |G_k(v_n)|}\le   \io\big (|H (x,\nabla v_n)|+|f|\big)(e^{\mu_n  |G_k(v_n)|}-1)\\
\le  p_n K_1\io |\nabla G_k(v_n) |^{p_n}e^{\mu_n  |G_k(v_n)|}+ \left(K_2+\|f\|_{\elle{\infty}}\right)\io(e^{\mu_n  |G_k(v_n)|}-1)\,, 
\end{split}
\]
that implies 
 \be\label{expest}
 \io|\nabla G_k(v_n)|^{p_n}   +  (k-A) \int_{\{|v_n|>k\}} (e^{\mu_n  |G_k(v_n)|}-1)\le 0. 
\ee
If we choose any $k>A$ we get a contradiction, since the left hand side above turns out to  be strictly positive.  
Consider now $k=0$ in \eqref{expest}, then  we  get
\[
\io|\nabla v_n|^{p_n} \le (K_2+\|f\|_{\elle{\infty}}) \io(e^{\mu_n  |v_n|}-1)\le B_n^{p_n}.
\]
Gathering together the above estimates, we  obtain \eqref{pestimate}. 
\end{proof}

Now we are ready to prove the main result of the section.
\begin{proof}[Proof of Theorem \ref{33}]
Thanks to   \eqref{pestimate}, we have that
\[
\|v_n\|_{\elle{\infty}}\le A \quad \mbox{and} \quad \|\nabla v_n \|_{\elle{p_n}}\le B_{n}.
\]
Using H\"older's inequality, it follows that for any $N<q<p_n$ we have
\[
\|\nabla v_n \|_{\elle{q}}\le  \|\nabla v_n \|_{\elle{p_n}}|\Omega|^{\frac1q-\frac 1{p_n}}.
\]
Since $B_{n}\to B_{\infty}$ (see Remark \ref{limitB}) and thanks to the $L^{\infty}$ bound, we deduce that the sequence $\{v_n\}$ is bounded in $W^{1,q}(\Omega)$ and, up to a (not relabeled) sequence, we get that $v_n\to v$ weakly in $W^{1,q}(\Omega)$ and uniformly in $\Omega$. The lower semicontinuity of the norm also implies that 
\[
\|\nabla v \|_{\elle{q}}\le  B_{\infty}\,, \qquad \forall q>1\,.
\]
Therefore we conclude that $v(x)$ satisfies 
$$\|v\|_{W^{1,\infty}(\Omega)}\le A+B_{\infty}.$$

Let us now focus on the equation solved by $v$. We start by proving that 
\be\label{subsol}
\max\big\{1-|\nabla v|,-\Delta_{\infty}v - H_{\infty}(x,\nabla v)     \big\}\le0  \qquad \mbox{ in } \Omega_v^+ \,.
\ee

Indeed, let us pick any $x_0\in  \Omega_v^+ $ and any $\vp \in C^{2} (\Omega)$ such that $v (x)-\vp(x)$ has a strict local maximum at $x_0$ with $v (x_0)-\vp(x_0)= 0$.  We need to check that 
\[
\max\big\{1-|\nabla \varphi(x_0)|,-\Delta_{\infty}\varphi(x_0)-H_{\infty}(x_0,\nabla\varphi (x_0)) \ \big\}\le0.
\]
Let us recall  that $v$ is the (uniform) limit of   solutions to \eqref{drift}, so that  there exist a sequence of real numbers $\eps_n\to0$ and a sequence of points $x_n\to x_0$ such that $v_n(x)-\vp(x)-\eps_n$ has a local maximum at $x_n$ and $v_n(x_n)-\varphi(x_n)=\eps_n$.

Since both $f $ and $H_{\infty}$ belong to $C^0 (\overline{\Omega})$, then (see Theorem 1 in \cite{Lieb})   any  $v_n$ belongs to $C^{1,\alpha} (\overline{\Omega})$ for some $\alpha \in (0,1)$ (see also \cite{Dib}) and thus $v_n$ turns out to be also a viscosity solution to \eqref{drift.intro}. Hence   we have that
\be\label{taglia}
\varphi(x_n)+\eps_n - f(x_n)- \Delta \vp (x_n)|\nabla \vp (x_n) |^{p_n-2} - (p_n-2) \Delta_{\infty} \vp (x_n) |\nabla \vp (x_n) |^{p_n-4}- H_{n}(x_n,\varphi (x_n)) \le 0\,.
\ee
Assume at first that $1-|\nabla\vp(x_0)| \leq 0$. It follows that $|\nabla\vp(x_n)|>0$ 
for  $n$ large enough and since $x_0$ belongs to the (open) set $ \Omega_v^+ $, then 
$\varphi(x_n)+\eps_n - f(x_n)$  is strictly positive and we can drop it from the previous inequality \eqref{taglia}. 
Hence, we divide \eqref{taglia} by  $(p_n-2)|\nabla \vp (x_n) |^{p_n-4}$ and we take the limit with respect to $n$: exploiting \eqref{Hconvergence} we deduce that 
\[
-\Delta_{\infty}\vp(x_0) - H_{\infty}(x_0,\nabla\varphi (x_0))  \le 0.
\]

Conversely,  if $|\nabla \vp(x_0)|<1$, we can directly pass  to the limit in \eqref{taglia}, and we get $\vp(x_0)-f(x_0)=u(x_0)-f(x_0)\le0$, that is in contradiction with the case under consideration.

Next we prove that 
\be\label{supersol}
\max\{1-|\nabla v|,-\Delta_{\infty}v -H_{\infty}(x,\nabla v)  \big\}\ge0 \qquad \mbox{ in } 
\{v(x) -f (x) >0 \}.
\ee
Consider  any $x_0\in \Omega_v^+$ and take any $\vp \in C^{2} (\Omega)$ such that $v (x)-\vp(x)$ has a strict local minimum at $x_0$ with $v (x_0)-\vp(x_0)=0$. We want to prove that for such a $\vp$ we have that 
\[
\max\big\{1-|\nabla \varphi(x_0)|,-\Delta_{\infty}\varphi(x_0)-H_{\infty}(x_0,\nabla\varphi (x_0)) \ \big\}\ge0.
\]
If $1-|\nabla\varphi(x_0)|\ge0$ the inequality is trivially satisfied, so we deal with  the case $|\nabla \varphi(x_0)|>1$. 
As before, we recall that  there exist a sequence of real numbers $\eps_n\to0$ and a sequence of points $x_n\to x_0$ such that $v_n(x)-\vp(x)-\eps_n$ as a local minimum at $x_n$ and $v_n(x_n)-\varphi(x_n)=\eps_n$ such that
\[
\vp (x_n)+\eps_n -f(x_n) - \Delta \vp (x_n)|\nabla \vp (x_n) |^{p_n-2} - (p_n-2) \Delta_{\infty} \vp (x_n) |\nabla \vp (x_n) |^{p_n-4} - H_{n}(x_n,\varphi (x_n)) \ge 0.
\]
Dividing the above inequality by $(p_n-2)|\nabla \vp (x_n) |^{p_n-4}$ and taking the limit with respect to $n$,  we conclude that
\[
-\Delta_{\infty}\vp(x_0) - H_{\infty}(x_0,\nabla\varphi (x_0))  \le 0.
\]

Gathering  \eqref{subsol} and \eqref{supersol}, we conclude that $v$ is a viscosity solution to
\[
\max\big\{1-|\nabla v|,-\Delta_{\infty}v - H_{\infty}(x,\nabla v)     \big\}=0  \qquad \mbox{ in }  \Omega_v^+ \,.
\]

The proof that $v$ solves 
\[
\min\big\{|\nabla v|-1,-\Delta_{\infty}v - H_{\infty}(x,\nabla v)     \big\}=0  \qquad \mbox{ in }  \Omega_v^-\,, 
\]
follows similarly to the previous case.  
 
 Let us focus now on the boundary conditions. For the sake of brevity we only prove the first one of the  four inequalities \eqref{bc1}--\eqref{bc2}, since the proofs of the other ones  follow in the same way.

Assume that there exists $\varphi\in C^2(\overline{\Omega})$ such that $v-\varphi$ has a strict local maximum at $x_0$, with  $v(x_0)=\varphi(x_0)$. 
Since $\{v_n\}$ converges uniformly in $\overline{\Omega}$  to $v$, there exist a sequence of real numbers $\eps_n\to0$ and a sequence of points $ \bar{\Omega} \ni x_n\to x_0$ such that $v_n-\vp-\eps_n$ attains a local maximum at $x_n$ and $v_n(x_n)=\varphi(x_n)+\eps_n$.

Let us assume that (up to a subsequence) $\{x_n\}\subset (\partial\Omega)^+_v$;   we recall that since any $v_n$ belongs to $C^{1,\alpha}(\overline{\Omega})$ then 
$$\frac{\partial v_n (x)}{\partial\nu} = 0,  \qquad\forall x \in \partial \Omega.$$ 
Using that $\varphi(x)+\eps_n$ touches $v_n$ from above at $x_n$, we have that
\[
\frac{\partial \varphi(x_n)}{\partial\nu}\le\frac{\partial v(x_n)}{\partial\nu}=0.
\]
Passing to the limit as $n\to\infty$, it follows 
\[
\frac{\partial \varphi(x_0)}{\partial\nu}\le 0.
\]
On the contrary, if we  assume that  $\{x_n\}\subset \Omega^+_v$, we take advantage to the fact that $v_n$ is a viscosity subsolution to \eqref{aux},  i.e. inequality \eqref{taglia} holds true. Going back to the proof of Theorem \ref{33} we observe that    $|\nabla \varphi(x_0)|<1$ yields to a a contradiction and that $|\nabla \varphi(x_0)|-1 \geq 0$ implies that   
$
-\Delta_{\infty}\varphi (x_0)- H_{\infty}\big(x,\nabla \varphi(x_0) \big)  \leq 0\,,$ as desired.
\end{proof}

\section{The limit equation for the model problem}\label{4}

Now, let us go back to the limit as $p\to \infty$ to \eqref{drift.intro} with $H\equiv0$. 

%Our previous results give that the unique minimizer is a solution to a double eikonal equation.
%
%\begin{teo}
%Assume $f\in C(\overline{\Omega})$ and let $u\in \mathbb{X}(\Omega)$ be the unique minimizer to \eqref{16:50} provided by Proposition \ref{auxlem}. Then it solves the following  equation
%\be\label{eqlim}
%\begin{cases}
% 1-|\nabla u| =0 &\mbox{in} \ \Omega^+_u \,,\\
%|\nabla u|-1 =0 &\mbox{in} \ \Omega^-_u \,,
%\end{cases}
%\ee
%and it satisfies the following boundary conditions
%\be\label{boundarycon}
%\begin{cases}
%1 - |\nabla  u |   \le 0 &\mbox{on} \ (\partial\Omega)^+_u \, , \\
% \displaystyle \max\left\{1 - |\nabla  u |,\frac{\partial u }{\partial \nu }\right\}\ge0 & \mbox{on } (\partial\Omega)^+_u \, ,
%\end{cases}
%\quad
%\mbox{and}
%\quad
%\begin{cases}
% |\nabla  u | -1  \ge 0 &\mbox{on} \ (\partial\Omega)^-_u \, , \\
% \displaystyle \min\left\{|\nabla  u | -1,\frac{\partial u }{\partial \nu }\right\}\le0 & \mbox{on } (\partial\Omega)^-_u \, .
%\end{cases}
%\ee
%\end{teo} 

\begin{proof}[Proof of Theorem \ref{16:15}]
Since $u\in \mathbb{X} (\Omega) $, by Proposition \ref{L1} we have  that both
\[
1-|\nabla u|\ge0 \qquad \mbox{and} \ \ \ |\nabla u|-1 \le 0 \qquad \mbox{in} \  \Omega.  
\]
The reverse inequalities in $\Omega^+_u$ and $\Omega^-_u$ follow directly from Theorem \ref{33}.
 
As far as the boundary conditions are concerned, notice that inequality 
\[
\max\left\{1 - |\nabla  u |,\frac{\partial u }{\partial \nu }\right\}\ge0 \quad \mbox{on } (\partial\Omega)^+_u \,,
\]
follows directly from Theorem \ref{33}. In order to prove that 
\be\label{2:23}
1 - |\nabla  u |   \le 0 \mbox{ on } \ (\partial\Omega)^+_u \,,
\ee
assume by contradiction that there exists $\varphi\in C^2$ that touches ${u}$ from above at $x_0$ and that verifies 
\[
1-|\nabla \varphi(x_0)|>0.
\]
Take a small $r$ such that $B_r(x_0)\cap\partial \Omega\subset (\partial\Omega)^+$ and $|\nabla \varphi  (y)|<1$ for all $|y-x_0|\le r$ and $B_r (x_0)$.
Therefore, set
$$
\mu  \vcentcolon = \dys \inf_{x \in \partial B_{r } (x_0) \cap \Omega} \big( \varphi  (x) - {u} (x) \Big)>0 \quad \mbox{and} \quad \eta \vcentcolon= \dys \inf_{x \in    B_{r } (x_0)\cap \Omega } \big( \varphi  (x) - f (x) \Big) >0
$$
and consider, for $\varepsilon\vcentcolon =\min\{\mu , \eta \}$, the following function
\[
\tilde u=\begin{cases}
\min\{ u,\varphi -\varepsilon \} &\mbox{in} \ B_r(x_0)\cap \overline\Omega \, , \\
 u&\mbox{in} \ \overline\Omega\setminus B_r(x_0)\,.
\end{cases}
\]  
Thanks to the definition of $\varepsilon$, $\tilde u$ is continuous, $\tilde u> f$ in $A^+$ and $\tilde u(x_0)<u(x_0)$. Moreover the choice of $r$ also implies that $\tilde u\in\mathbb{X} (\Omega) $. This would imply $J_{\infty}(\tilde u)<J_{\infty}(u)$, that contradicts the minimality of $u$. Therefore \eqref{2:23} holds true.

To deal with the boundary condition on $(\partial\Omega)^-$, we can follow exactly the same strategy and hence we omit the details. 
\end{proof}

Even if the minimizer of \eqref{16:50} is unique,  in general  problem 
\eqref{eqlim.intro}--\eqref{boundarycon.intro} possesses more than one solution.  
In the following we   provide a solution to \eqref{eqlim.intro}, that satisfies the same boundary condition, and  that differs from the minimizer of $J_{\infty}$ (unless $f \in  \mathbb{X} (\Omega)$).

\begin{proof}[Proof of Theorem \ref{obstacle}]
Let us consider the set
\[
M=\{v\in \mathbb{X} (\Omega) \ :  \  \ v\ge  f    \},
\]
and notice that it is not empty since the constant $\|f\|_{L^{\infty} (\Omega)} \in  M$, and let us set
\[
\alpha=\inf \left\{\io v \ : \ v\in M\right\}.
\]
By definition of $M$, we have that $\alpha\geq \io f$ and we   consider a   sequence $\{v_n\}\subset M$ whose integral on $\Omega$ converges to $\alpha$. Observe that for any $n\in\mathbb{N}$, we have that $v_n\ge f$ and moreover one can assume, without loss of generality, that $v_n\le \|f\|_{L^{\infty} (\Omega)}$ (otherwise it is sufficient to truncate  $v_n$ at high  $ \|f\|_{L^{\infty} (\Omega)}$). Hence the sequence    $\{v_n\}$ is equibounded. 

Moreover, since  $\{ v_n\} \subset \mathbb{X} (\Omega)$, then it is  also  equicontinuous. Thus by   Ascoli Arzel\`a Theorem, we have that 
$ v_n  $ converges  uniformly in $\overline{\Omega}$ to some $\overline u \in \mathbb{X}$,   that satisfies $\alpha=  \io \overline u $.

Next,  we  prove that $\overline u$ turns out to be  a viscosity solution to \eqref{super}. Thanks to Proposition \ref{L1}, we immediately deduce that $  \overline u$ satisfies 
\[
 1-|\nabla \overline u|  \geq 0 \ \ \  \mbox{in } \Omega .
\]
Hence, to conclude the proof, we take any  $x_0\in  \overline{\Omega^+_u} \cap(\partial\Omega)^+_u $ and assume by contradiction that there exists $\varphi\in C^2$ that touches $\overline{u}$ from above at $x_0$ and that satisfies 
\[
1-|\nabla \varphi(x_0)|>0.
\]
Since $\varphi$ is smooth, there exists a small $r>0$ such that  $|\nabla \varphi (y)|<1$ for all $y\in B_r (x_0) \cap \overline{\Omega}$.
Therefore, we  set
$$
\mu \vcentcolon = \dys \inf_{x \in \partial B_{r } (x_0)\cap \Omega} \big( \varphi  (x) - \overline{u} (x) \Big)>0 \qquad \mbox{and} \qquad 
\eta \vcentcolon= \dys \inf_{x \in  B_{r } (x_0)\cap \Omega} \big( \varphi  (x) - f (x) \Big) >0
$$
and we consider
\[
\tilde u=\begin{cases}
\min\{\overline u,\varphi -\varepsilon \} &\mbox{in} \ B_r(x_0)\,,\\
 \bar u&\mbox{in} \ \Omega\setminus B_r(x_0)\,,
\end{cases}
\qquad \mbox{ with } \qquad \varepsilon\vcentcolon =\min\{\mu, \eta \}\,.
\]  
Thus,  thanks to the definition of $\varepsilon$, $\tilde u$ is continuous, $\tilde u> f$ in $\Omega^+_{\overline u}$ and $\tilde u(x_0)<\overline u(x_0)$. Moreover the choice of $r$ also implies that $\tilde u\in\mathbb{X} (\Omega) $.
However, by construction this would imply that $\io \tilde u<\alpha$, that yields to a contradiction.  Hence, we obtain that $1-|\nabla \varphi(x_0)|\leq 0$.

It remains to verify that, for any $x_0\in   \partial \Omega^+_{\overline u}$ and any $\varphi\in C^2$ that touches $\overline u$ from below at $x_0$, we have
\[
\max \left\{1 - |\nabla \varphi(x_0) |,\frac{\partial \varphi(x_0) }{\partial \nu }\right\}\ge0\,.
\]
Assume by contradiction that  $1 - |\nabla \varphi(x_0) |<0$ and $\frac{\partial \varphi(x_0) }{\partial \nu }=-\nabla \varphi(x_0) \cdot  \nabla \delta(x_0)<0$.

This latter condition implies that   $\nabla \varphi(x_0)$ directs  inside $\Omega$; moreover since $|\nabla \varphi(x_0) |>1$, we deduce that there exists $y_0\in\Omega$ (close enough to $x_0$) such that $\varphi(y_0)-\varphi(x_0)>|y_0 - x_0|$. Using that $\varphi$ touches $\overline u$ from below at $x_0$ we have that
\[
\overline u(y_0)-\overline u(y_0)\ge\varphi(y_0)-\varphi(y_0)>|y_0 - x_0|\,,
\]
that yields to a contradiction.

Finally, the uniqueness of the solution of the obstacle problem follows from the minimality of $\overline{u}$.
\end{proof}

\begin{remark}  
Analogously to \eqref{super} one can consider the solution to 
\be\label{sub}
\begin{cases}
\underline{u} \leq  f \qquad &\mbox{in } \Omega,\\
|\D \underline{u} (x)|-1 \leq 0 &\mbox{in } \Omega,\\  
 |\D \underline{u} (x)|-1 =0 &\mbox{if} \ \Omega^-_{\underline u} ,\\
 \displaystyle \min \left\{ |\nabla \underline{u} |-1,\frac{\partial \underline{ u } }{\partial \nu }\right\}\le0 & \mbox{on } (\partial\Omega)^-_{\underline u}\,.
\end{cases}
\ee
Arguing as before, we have that  also $\underline{u}$  solves \eqref{eqlim.intro}. 

\end{remark}

\begin{remark} {\rm
As already observed, we have that  in general 
$$
\underline{u} (x) \leq  f  (x) \leq \bar{u} (x) \qquad \mbox{ in } \Omega\,,
$$
while the variational solution $u$ is such that $u-f$ cannot have constant sign (unless $f\in \mathbb{X} (\Omega)$). 
}
\end{remark}

\begin{proof}[Proof of Proposition \ref{ifandonlyif}]
Let us define 
\be\label{cones}
\begin{matrix}
\medskip
\dys v^+(x):=\max_{y\in\partial \Omega^+_{w}} \big[f(y)-d(x,y) \big]\,, \quad  \forall   x\in \Omega^+_{w}\,, \\
\dys v^-(x):=\min_{y\in\partial \Omega^-_{w}} \big[f(y)+d(x,y)\big]\,,\quad \forall   x\in \Omega^-_{w}.
\end{matrix}
\ee

We prove, at first,  that $v^+$ solves $1-|\nabla v^+|=0$ and that $-\Delta_{\infty}v^+\le0$ in $\Omega^+_{w}$. Thanks to Lemma \ref{minmaxcone} (we recall that $v^+\in \mathbb{X} (\Omega^+)$) we already now that $1-|\nabla v^+|\ge 0$. Notice now that for any $y\in\partial \Omega^+_{w}$, the function 
$$
x\to f(y)-d(x,y) 
$$
is locally a cone (recall that the geodesic distance coincides locally with the euclidean one). Therefore, it solves the eikonal equation and it is infinity harmonic outside its vertex. Since  $v^+$ is defined as maximum of such functions, the stability properties of viscosity solutions (see Proposition 4.3 of \cite{CIL}) imply that  $1-|\nabla v^+|\le 0$ and $-\Delta_{\infty} v^{+}\le0$ in $ \Omega^+_{w}$.

A similar argument proves that $|\nabla v^-|-1=0$ and $-\Delta_{\infty} v^-\ge0$ in $ \Omega^-_{w}$. 

Consequently,  we found  that
\[
\begin{matrix}
\max\big\{1-|\nabla v^+|,-\Delta_{\infty}v^+\big\}=0, \ \ \ \mbox{and} \ \ \ 1-|\nabla v^+|=0 &\mbox{in} \  \Omega^+_{w} \,,\\[8pt]
\min \big\{|\nabla v^-|-1,-\Delta_{\infty}v^- \big\}=0 \ \ \ \mbox{and} \ \ \ |\nabla v^-|-1=0 &\mbox{in} \  \Omega^-_{w}.
\end{matrix}
\]
 
Let us assume now that $w$ solves
\begin{equation} \label{rrrr}
\begin{cases}
\max\big\{1-|\nabla w|,-\Delta_{\infty}w \big\}=0 \qquad &\mbox{in} \quad  \Omega^+_w \,,\\[8pt]
\min\big\{|\nabla w|-1,-\Delta_{\infty} w  \big\}=0 \qquad &\mbox{in} \quad  \Omega^-_w\,\,.
\end{cases}
\end{equation}
Thanks to the comparison principle (see Theorem 2.1 of \cite{J}) we have that 
$$w\equiv v^{\pm} \quad \mbox{  in } \ \ \Omega^{\pm}_w\,,$$ respectively. Thanks to the properties of $v^{\pm}$, we deduce that $w$ is also a solution to
\be\label{rree}
\begin{cases}
 1-|\nabla w| =0 &\mbox{in} \ \Omega^+_w\, ,  \\[8pt]
|\nabla w|-1 =0 &\mbox{in} \ \Omega^-_w \,.
\end{cases}
\ee
On the other hand if $w$ is a solution to \eqref{rree} then again for the comparison principle (see Theorem 5.9 of \cite{BC}) $w\equiv v^{\pm}$ in $\Omega^{\pm}_w$ and then it also solves \eqref{rrrr}.
\end{proof}

\section{A functional allowing jumps}\label{mumford}
We start this section recalling some basic facts abound the \emph{space of special functions of
bounded variation} $SBV(\Omega)$. For more details we address the interested reader to the overview \cite{fusco} or the classical monograph \cite{AFP}.\\
Let us recall that a function $v\in\elle 1$ belongs to $BV(\Omega)$ (bounded variation) if and only if there exists a vector valued Radon measure $Du$ such that
\[
\io v \, \mbox{div}(\varphi )=-\io \varphi  \, d (Dv) \ \ \ \forall \ \varphi \in C^{\infty}_c(\Omega,\mathbb{R}^N).
\]
Thanks to the Radon-Nikodym Theorem, we can uniquely decompose $Dv$ as $Dv=Dv^a+Dv^s$ with
\[
Dv^a\ll \mathcal{H}^N \quad \mbox{and} \quad Dv^s\perp\mathcal{H}^N.
\]
To give a better description of the absolutely continuous part and the singular one, we need to introduce some more concepts. For any $x\in\Omega$, we say that $v\in BV(\Omega)$ admits approximate limit $\bar{v}(x)\in\mathbb{R}$ if
\be\label{approxcont}
\lim_{\rho\to 0} \frac1{ |B_{\rho}(x)|}\int_{B_{\rho}(x)}|v(y)-\bar{v}(x)|dy=0
\ee
and approximate gradient $\nabla v(x)\in\mathbb{R}^N$ if 
\[
\lim_{\rho\to 0} \frac1{\rho |B_{\rho}(x)|}\int_{B_{\rho}(x)}|v(y)-\bar{v}(x)-\nabla v(x)\cdot(y-x)|dy=0.
\]
It can be proved (see Theorem 3.83 of \cite{AFP}) that any $v\in BV(\Omega)$ is approximately differentiable $a.e.$ and that $\nabla v\in (\elle1)^N$ is the density of the absolutely continuous part of $Dv$ with respect to the Lebesgue measure, namely
\[
Dv^a=Dv\llcorner \mathcal{H}^N=\nabla v\mathcal{H}^N.
\]
In order to describe also the singular part of $Dv$, let us consider the set $C_v\subset\Omega$ of all points where \eqref{approxcont} holds true and denote $S_u=\Omega\setminus C_v$. The set $S_v$ is called the set of approximated discontinuity of $v$ and it is $(N-1)$-countably rectifiable. Now, we introduce
\[
Dv^j=Dv\llcorner S_u \quad \mbox{and} \quad Dv^c=Dv\llcorner (\Omega\setminus S_u).
\]
The measure $Dv^j$ takes into account the jumps of the function $v$, while $D^c$ is the Cantor part of $Dv$. We notice that in \cite{AFP} $Dv^j$ is defined through the slightly smaller set $J_v$, the set of approximated jumps point of $v$; since $\mathcal{H}^{N-1}(S_v\setminus J_v)=0$ we keep using $S_v$ for the sake of simplicity (see also \cite{fusco}).
At this point it is easy to define the space of special functions of bounded variation as
\[
SBV(\Omega)=\{u\in BV(\Omega) \ : \ D^cu\equiv 0\}.
\]  
Now we are ready to define our functional $ \mathcal{I} :SBV(\Omega)\to \mathbb{R}$ as
\be\label{ms}
\mathcal{I}_p(v)=\frac1p\io|\nabla v|^p +\frac12\io|v-f|^2+\mathcal{H}^{N-1}(S_v).
\ee
Comparing \eqref{ms} with \eqref{mumsha}, we note that the new funcional framework allows to pass from a two variable functions to a single variable one, where the discontinuity set $K$ is replaced by the singular (jump) set $S_u$.

We also stress that $BV(\Omega)$ would not be a good ambient space for \eqref{ms} since the functional is not coercive with respect to the cantor part of gradient measure.
To be convinced that $SBV(\Omega)$ is the right ambient space to settle the minimization of $\mathcal{I} _p$, it is enough to look at the following result.

\begin{teo}{[Theorems 4.7 and 4.8 of \cite{AFP}]}\label{compact}
Let $\{\un\}\subset SBV(\Omega)$ such that
\[
\io|\nabla v_n|^{p}+\mathcal{H}^{N-1}( S_{v_n})+\|v_n\|_{\elle{\infty}}\le C.
\]
Then there exists $v\in SBV(\Omega)$ such that, up to a subsequence $v_n\to v$ strongly in $\elle 1$, $\nabla v_n\rightharpoonup \nabla v$ weakly in $\elle p$, $D v_n\rightharpoonup D v$ weakly star in the sense of measure and
\[
\io |\nabla v|^p\le \liminf_{n\to\infty}\io|\nabla v_n|^{\pn}
 \ \ \ \mbox{and} \ \ \ 
\mathcal{H}^{N-1}( S_{v})\le \liminf_{n\to\infty} \mathcal{H}^{N-1}( S_{v_n}).
\]

\end{teo}
Theorem \ref{compact} easily implies the existence of a minimizer $u_p\in SBV(\Omega)$ for \eqref{ms}. Indeed it is easy to see that any minimizing sequence $v_n$ can be chosen so that $\|v_n\|_{\elle{\infty}}\le \|f\|_{\elle{\infty}}$, basically because truncations with a large constant always decrease the energy of the functional.

What is not clear is that if the couple $(S_{u_p},u_p)$ is a minimizer of $\eqref{mumsha}$. The main difficulty here is that in general the singular set of a $SBV$ function need not to be closed and can be even dense in $\Omega$.
What is possible to prove is that any minimizer $u_p\in SBV(\Omega)$ of \eqref{ms} satisfies $\mathcal{H}^{N-1}(\overline{S}_{u_p}\setminus S_{u_p})=0$ and that indeed $(\overline{S}_{u_p},u_p)$ is a minimizer of $\eqref{mumsha}$. We gather all this information in the following result.

\begin{teo}{[Theorems 7.15, 7.21 and 7.22 of \cite{AFP}]}
Set $p\in(1,\infty)$, $r\in[1,\infty)$ and $f\in \elle{\infty}$. Then there exists a minimizer $u_p\in SBV(\Omega)$ to \eqref{ms}. Moreover there exists $\theta$ depending on $N$ and $p$ such that
\[
\mathcal{H}^{N-1}(S_{u_p}\cap B_{\rho}(x))>\theta \rho^{N-1}
\]
for any $x\in \overline{S_{u_p}}$ and $\rho>0$.

This lower bounds implies that $\mathcal{H}^{N-1}(\overline{S_u}\setminus S_u)=0$. Thus the pair given by $K=\overline{S}_{u_p}$ and $u\in W^{1,p}(\Omega\setminus \overline{S}_{u_p}) $ and is a minimizer of the \emph{classical} Mumford–Shah problem \eqref{mumsha}.
\end{teo}

Our next step is to prove that any sequence $\{u_n\}\subset SBV(\Omega)$ of minimizers of \eqref{ms}, associated to a sequence $\pn\to \infty$, converges (up to a subsequence) to a minimizer of 
\be\label{msinf}
M_{\infty}(v)=\frac1r\io|v-f|^r+\mathcal{H}^{N-1}(S_v) \ \ \ \mbox{with} \ v\in SBV(\Omega), \ |\nabla v|\le1.
\ee

\begin{teo}\label{Gammacon}
Let $\{u_n\}\subset SBV(\Omega)$ be the sequence of minimizer of the functional \eqref{ms} relative to a diverging sequence of real number $p_n$. Then any $u_{\infty}\in SBV(\Omega)$, accumulation point of $\un\to u$ (in $\elle 1$), is a minimizer of \eqref{msinf}.
\end{teo}
\begin{proof}
Recalling that $M_{\pn}(\un)\le M_{\pn}(0)$ and that $\|\un\|_{\elle{\infty}}\le \|f\|_{\elle{\infty}}$, it results
\[
\frac1{\pn}\io|\nabla \un|^{\pn} +\|\un\|_{\elle{\infty}}+\mathcal{H}^{N-1}(S_{\un}) \le C\|f\|_{\elle{\infty}}.
\]
Using H\"older inequality we get
\[
\left(\io|\nabla \un|^q\right)^{\frac1q}\le \left(p C\|f\|_{\elle{\infty}}\right)^{\frac1{\pn}}|\Omega|^{\frac1q+\frac1{\pn}}
\]
for any $q<p_n$.
Theorem \ref{compact} and a diagonal argument then implies that there exists $u_{\infty}\in SBV(\Omega)$, with $|\nabla u_{\infty}|\le 1$ such that $\un\to u_{\infty}$ strongly in any $\elle q$, $\nabla \un\rightharpoonup \nabla u_{\infty}$ weakly in any $(\elle q)^N$, with $q\in[1,\infty)$,  $D \un\rightharpoonup D u_{\infty}$ weakly star in the sense of measure and
\be\label{cip}
M_{\infty}(u_{\infty})\le \liminf_{n\to \infty} M_{\pn}(\un).
\ee
Now we claim that 
\[
M_{\infty}(u_{\infty})\le M_{\infty}(z) \ \ \ \forall \ z\in SBV(\Omega)\cap\elle{\infty}, \ |\nabla z|\le 1.
\]
This concludes the proof since, by a simple truncation argument, if reached, the infimum of $J_{\infty}$ in $SBV(\Omega)\cap\elle{\infty}$ has to be bounded. By definition of $u_n$ we have that
\be\label{ciop}
M_{p_n}(\un)\le M_{p_n}(z )\le \frac1{\pn}|\Omega|+\frac12\io|z-f|^2+\mathcal{H}^{N-1}(S_{z}).
\ee
Taking the liminf with respect to $n$  we obtain the desired result.
\end{proof}

Notice that in the proof of Theorem \ref{Gammacon} we have proved that the functional $M_{\pn}$ $\Gamma$-converges to the functional $M_{\infty}$. 
For completeness let us recall the definition of $\Gamma$-convergence.

\begin{defin}\label{gaconv}
Let $F_n\,: \mathbb{X} \to \mathbb{R}$ be a sequence of functionals. We say that $F_n$ $\Gamma$--converges to $F\,: \mathbb{X} \to \mathbb{R}$ if 
\begin{itemize}
\item[i)] for any $x\in \mathbb{X}$ and any $x_n \to x$, we have that 
$$
\liminf_{n \to +\infty} F_n (x_n) \geq F (x) \,; 
$$
\item[ii)] for any $x\in \mathbb{X}$,  there exists $x_n \subset \mathbb{X}$ such that  $x_n \to x$, and
$$
\limsup_{n \to +\infty} F_n (x_n) \leq F (x) \,. 
$$
\end{itemize}
\end{defin}
In our case $i)$ has been proved in \eqref{cip} while $ii)$ is trivially obtained taking the limit in the right hand side of \eqref{ciop}.

\section{Examples }\label{exa}
In this section we collect some explicit examples that clarify the behavior of the minimizers of some of the functional presendet in this manuscript.

\begin{examples} 
Suppose that 
$\Omega = (-1,1)$,   we are interested in considering the following minimization problem: 
$$
\inf_{\mathbb{X}(-1,1) }  J_{\infty}(v) \qquad  \mbox{with} \qquad \ J_{\infty}(v)=\frac12 \int_{-1}^1 (v-f)^2, 
$$
\begin{itemize} 

\item[Case 1. ]  Choose $f_1 (x) = k \chi_{(-r,r)} (x)$ with $r \in (0,1]$ and $k\in \R$ (see Figure 1).

We first observe that if $r=1$ then $f_1 \in \mathbb{X}(-1,1)$ and thus the choice $u_1 (x)  \equiv f_1 (x)$ is allowed. 

Conversely, if $r\in (0,1)$, then  $f_1 \not \in \mathbb{X}(-1,1)$, and we look for a solution $u_1 (x)$ of the form  
 $u_h (x) = \min\{ (h-|x|)_+, k\}$, with $h$ to be determined. 
 Thus $J_\infty(u_h)$ can be explicitly computed and it turns out that 
$$
J_\infty (u_h) = \int_{h-k}^r \Big( h- k -x \Big)^2 dx+ \int_r^{h} \Big( h-x \Big)^2 dx\,. 
$$
Minimizing its value with respect to $h$,  we have that the minimum is achieved at $h^*= r + \frac{k}{2}$. 

The heuristic that motivates this example is that in a neighbor  where $f$ is not smooth
the solution is a line with maximum slope allowed (i.e. $\pm1$) and far away it tries to be as close as possible to $f$. 

\bigskip

\item[Case 2. ]  Consider $f_2 (x) = 2|x| $ (see Figure 2).

In this case the picture is much more clear. Since $f_2\not\in \mathbb{X}(-1,1)$,  the closest function to $f_2$ in $\mathbb{X}(-1,1)$  is a line with slope $+1$ in $(0,1]$ and with slope $-1$ in $[-1,0]$. 
Since $f_2$ is even, $u_2$ inherits the same property, so that 
$$
J_\infty (u_h) = \int_{0}^1 \Big(2x-(x+h) \Big)^2 dx \,,
$$
attains its minimum value for  $h=\frac12$. 

\bigskip

\item[Case 3. ]  Consider $f_3(x)=\sqrt{|x|}$ (see Figure 3).

By symmetry it follows  that the solution is even and that it has the   form
\[
u_3 (x)=
\begin{cases}
|x|-s+\sqrt{s} & \mbox{for} \ |x|\leq s \\
\sqrt{|x|}& \mbox{for} \ |x| > s \,,
\end{cases}
\]
with $s>\frac14$. To determine the value of $s$ we have to compute
\[
\begin{split}
\frac{\partial }{\partial s}J_{\infty}(u)&=\frac{\partial }{\partial s}\int_0^{s}(\sqrt{|x|}-|x|+s-\sqrt{s})^2dx\\
&=2\left(1-\frac{1}{2\sqrt{s}}\right)\int_0^{s}(\sqrt{|x|}-|x|+s-\sqrt{s})dx=2\left(1-\frac{1}{2\sqrt{s}}\right)\left(\frac{s^2}{2}-\frac{s^{\frac32}}{3}\right).
\end{split}
\]
Therefore we deduce that the minimum is achieved for $s=\frac{4}{9}$.
\end{itemize}
\end{examples}

\begin{example}
Let us underline the importance of minimizing $J_\infty $  in $SBV (\Omega)$ by showing that for certain $f$, when it is allowed, the minimizer does jump. 

In order to avoid technicalities, we present first full details in dimension 1, being the extension of these ideas to any dimension quite similar.

As already observed in the previous example the minimizer of 
$$
 J_\infty (u)= \frac12 \int_{-R}^R \Big(u(x)- k\chi_{ (-r,r)} \Big)^2 dx \qquad \mbox{    with    } u \in \mathbb{X} \big(-R,R\big)  
$$
is  $$u_h (x) = \min\big\{ (h-|x|)_+, k\big\},$$ 
with $h^*= r + \frac{k}{2}$. In addition, we have 
$$J_\infty (u_{h^*}) = \frac1{12} k^3.$$

Observe that, since such an $f \in SBV (-R,R)$,  we can consider 
$$
\mathcal{I} (v) = \frac12 \int_{-R}^R \Big(u(x)- k\chi_{ (-r,r)} \Big)^2 dx + \mathcal{H}  (S_u) \qquad \mbox{  with } u \in SBV \big(-R,R\big)   
$$
with $ \mathcal{H}  (S_u) $ the measure of the singular set of $u$. 
Observe that  $\mathcal{H}  (S_f)=2 $ and since $f \in SBV (-R,R)$ thus $\mathcal{I} (f) = 2 $. Consequently if 
$$2 < \frac1{12} k^3,$$ we have 
$$\mathcal{I}(f) < J_\infty (u_h^*) ,$$
and we conclude that the minimizer indeed has jumps.

 As far as the case in dimension $N\geq 2$ is concerned, we repeat  the same argument: 
 in the domian    
 $\Omega = B_R (0)$ we set   $$f(x) = k \chi_{B_r (0)}(x)$$ with   $r>0$ to be chosen large enough. 
 
Thus the minimizer of 
$$
 J_\infty (u)= \frac12 \int_{B_R (0)} (u(x)- k\chi_{B_r (0)} )^2 dx \qquad \mbox{    with    } u \in \mathbb{X} \big(B_R (0)\big)  
$$
has the form
$$
u_h (x) = \min \big\{k,(h -|x|)_+ \big\}\,.
$$

The minimum value of $J_\infty (u_h)$ is achieved for a suitable value $h^*>0$ that is not so easy to be computed explicitly and which 
satisfies 
$$\dys \lim_{r \to \infty} \frac{r}{h^*} = 1.$$
Thus, tedious but not difficult computations yield to 
$$
J_\infty (u_{h^*}) \sim \frac{k^3}{24} \omega_N  r^3 \qquad \mbox{ as } r \to +\infty\,.
$$
Analogously, it is not hard to see that  for the functional
$$
\mathcal{I}  (u) =  \frac12 \int_{B_R (0)} (u(x)- k\chi_{B_r (0)} )^2 dx  + \mathcal{H}^{N-1} (S_u) \qquad \mbox{  with } u \in SBV (\Omega) 
$$
since $f \in SBV (\Omega) $ we have that 
$$
\mathcal{I}  (f)= N\omega_N r^{N-1}\,.
$$

Consequently, if $r$ is large enough and $k^3 >24  N$, then it holds that  
$$
\mathcal{I}  (f) = N\omega_N r^{N-1} <  J_\infty (u_0),
$$
and the minimizer jumps.
\end{example}


\begin{thebibliography}{999}

\bibitem{ALT} A. Alvino, P. L. Lions, G. Trombetti. {\it Comparison results for elliptic
and parabolic equations via Schwarz symmetrization}, Ann. Inst. H. Poincaré
Anal. Non Linéaire, {7} (1990) 37–65.

\bibitem{AFP} L. Ambrosio, N. Fusco, D. Pallara. 
{\it Functions of bounded variation and free discontinuity problems},
 Oxford Mathematical Monographs, The Clarendon Press Oxford University Press, New York, 2000.

\bibitem{ar3}   G. Aronsson.  
{\it Extension of functions satisfying Lipschitz conditions}, Ark. Mat. 6 (1967), 551--561.

\bibitem{acj}  G. Aronsson, M. G. Crandall, P. Juutinen.
{\it A tour of the theory of absolutely minimizing functions},
Bull. Amer. Math. Soc. (N.S.) 41 (2004), no. 4, 439--505.

\bibitem{BC} M. Bardi, I. Capuzzo-Dolcetta. 
{\it Optimal Control and Viscosity Solutions of Hamilton-Jacobi-Bellman Equations}, 
Modern Birkh\"auser Classic, 1997.

\bibitem{bej} E. N. Barron, L. C. Evans, R. R. Jensen.
{\it The infinity Laplacian, Aronsson's equation and their generalizations},
Trans. Amer. Math. Soc. 360 (2008), 77--101.

\bibitem{bjw1}  E. N. Barron, R. R. Jensen, C. Y. Wang.
{\it The Euler equation and absolute minimizers of $L\sp \infty$ functionals}, 
Arch. Ration. Mech. Anal. 157 (2001),  255--283.

\bibitem{bfm} M.F. Betta, V. Ferone, A. Mercaldo. {\it Regularity for solutions of nonlinear elliptic
equations}, Bull. Sci. Math. 118  (1994), 539--567.

\bibitem{BBM}   T. Bhattacharya, E. Di Benedetto, J. Manfredi. 
{\it Limits as $p \to \infty$ of $\Delta_p u_p = f$ and related extremal problems},  
Rend. Sem. Mat. Univ. Politec. Torino, (1991), 15--68.



\bibitem{lucio} L. Boccardo. {\it Some developments on Dirichlet problems with discontinuous coefficients}, Boll. Unione Mat. Ital. 9,  (2009), 285–297.

\bibitem{nsesamai} L. Boccardo, S. Buccheri, G.R. Cirmi, {\it Calderon–Zygmund–Stampacchia theory for infinite energy solutions of nonlinear elliptic equations with singular drift}, Non. Diff. Equations Appl. NoDEA, 27 (2020), no. 4, Paper No. 38, 17 pp.

\bibitem{BMP1} L. Boccardo, F. Murat, J.-P. Puel. 
{\it  R\'esultats d’existence pour certains probl\`emes elliptiques quasilin\'eaires,}  Annali della Scuola Normale Superiore di Pisa, Classe di Scienze, 11 (1984), 213--235.

\bibitem{BMP4} L. Boccardo, F. Murat, J.-P. Puel.  {\it Quelques propri\'et\'es des op\'erateurs elliptiques quasi lin\'eaires,}, C. R. Acad. Sci. Paris. {307}, 1988.

\bibitem{BMP} L. Boccardo, F. Murat, J.-P. Puel. 
{\it $L^\infty$ estimate for some nonlinear elliptic partial differential equations and application to an existence result}, SIAM J. Math. Anal. {23}, (1992), 326--333.


%\bibitem{BCP} A. Braides, V. Chiad\'o Piat.  {\it Integral representation results for functionals defined on $SBV(\Omega, \mathbb{R}^m)$}, J. Math. Pures Appl. 75 (1996), 595--626.

 \bibitem{Br} H. Brezis.  
 {\it  Functional analysis, Sobolev spaces and partial differential equations},  Universitext. Springer, New York, 2011. xiv+599


 \bibitem{abilitazione} S. Buccheri. {\it The Bottaro-Marina slice method for distributional solutions to elliptic equations with drift term},  AIP Conference Proceedings 2425, 210001  (2022). 

\bibitem{BL} S. Buccheri, T. Leonori. 
{\it Large solutions to quasilinear problems involving the $p$-Laplacian as $p$ diverges},  Calc. Var. Partial Differential Equations, 60 (2021), no. 1, Paper No. 30, 23 pp. 


\bibitem{BLR} S. Buccheri, T. Leonori, J. D. Rossi. 
{\it Strong convergence of the gradients for $p$-Laplacian problems as $p\to \infty$}, J. Math. Anal. Appl. 495 (2021),   Paper No. 124724, 11 pp. 

\bibitem{DCL} E. De Giorgi, M. Carriero, A. Leaci. 
{\it  Existence theorem for a minimum problem with free discontinuity set}, 
Arch. Rational Mech. Anal. 108 (1989), 195--218.

\bibitem{Dib} E. DiBenedetto. 
{\it $C^{1+\alpha}$ local regularity of weak solutions of degenerate elliptic equations}, Nonlinear Anal. 7,  (1983), 827--850.
 
 \bibitem{cdp}   T. Champion,  L. De Pascale.
{\it A principle of comparison with distance functions for absolute minimizers}, 
J. Convex Anal.,  14,  (2007),  515--541.

\bibitem{cra} M. G. Crandall.  
{\it An efficient derivation of the Aronsson equation},
Arch. Ration. Mech. Anal. 167, (2003), 271--279.

\bibitem{ceg}  M. G. Crandall, L. C. Evans,   R. F. Gariepy. 
{\it Optimal Lipschitz extensions  and the infinity Laplacian}, 
Calc. Var. Partial Differential Equations, 13 (2001), 123--139.

\bibitem{CIL} M. Crandall, H. Ishii, P.L. Lions. 
{\it User's guide to viscosity solutions of second order partial differential equations},  
Bulletin Amer. Math. Soc., {27}, (1992), 1--67. 

\bibitem{EG}   L.C. Evans,  W. Gangbo. 
{\it Differential equations methods for the Monge-Kantorovich mass transfer problem}, 
Mem. Amer. Math. Soc., 137 (1999), no. 653.

\bibitem{fusco} N. Fusco. 
{\it An overview of the Mumford-Shah problem}, 
Milan J. Math. 71, (2003) 95--119.

\bibitem{J} R. Jensen. 
{\it Uniqueness of Lipschitz extension: minimizing the sup norm of the gradient}, 
Arch. Rat. Mech. Anal. 123, (1993) 51--74.

\bibitem{LPR} T. Leonori, A. Porretta, G. Riey. 
{\it  Comparison principles for p-Laplace equations with lower order terms},  
Ann. Mat. Pura Appl. 196, (2017),  877--903. 

%\bibitem{LL} J. Leray, J.L. Lions. {\it  Quelques r\'esultats de Visik sur les probl\'emes elliptiques non lin\'eaires par les m\'ethodes de Minty-Browder}, Bull. Soc. Math. France, 93 (1965), 97--107.

\bibitem{Lieb} G. Lieberman. 
{\it Boundary regularity for solutions of degenerate elliptic equations}, 
Nonlinear Anal. 12, (1988),  1203--1219.

\bibitem{MO}  M. Medina, P. Ochoa. 
{\it On viscosity and weak solutions for non-homogeneous $p$-Laplace equations},   Adv. Nonlinear Anal. 8, (2019), 468--481.

\bibitem{MN} T. Mengesha, N. C. Phuc. {\it Quasilinear Riccati type equations with distributional data in Morrey space framework}, J. Differential Equations, 260, (2016) 5421--5449.

\bibitem{rs} J. D. Rossi,  N. Saintier.   
{\it On the first nontrivial eigenvalue of the $\infty$--Laplacian with Neumann boundary conditions},  Houston J. Math. 42, (2016),  613--635.

\end{thebibliography}
\end{document}